\documentclass[twoside]{article}

\usepackage[accepted]{aistats2021}
\usepackage{graphicx}
\usepackage{subfigure}
\usepackage{booktabs} 
\usepackage{amssymb}
\usepackage{enumitem}
\usepackage{epstopdf,algorithm,algorithmic}
\usepackage{color}
\usepackage[algo2e,linesnumbered,vlined,ruled]{algorithm2e}
\usepackage[round,sort]{natbib}

\usepackage{hyperref}
\input{mysymbol.sty}
\newtheorem{lemma}{\bf  Lemma}
\newtheorem{theorem}{\bf  Theorem}
\newtheorem{assumption}{\bf Assumption}

\newtheorem{proposition}{\bf  Proposition}

\usepackage[ruled,vlined]{algorithm2e}
\begin{document}

%
\runningtitle{Convergence Rates of  Average-Reward MARL}

%

\twocolumn[

\aistatstitle{Convergence Rates of  Average-Reward Multi-agent \\Reinforcement Learning via Randomized Linear Programming}

\aistatsauthor{ Alec Koppel$^{*\dagger}$	 \And Amrit Singh Bedi$^{*\$}$ \And  Bhargav Ganguly$^\ddagger$  \And  Vaneet Aggarwal$^\ddagger$ }

\aistatsaddress{ $^{\dagger}$Amazon  \And $^{\$}$CISD, US Army Research Lab. \And $^\ddagger$Purdue University } ]

\begin{abstract}
In tabular multi-agent reinforcement learning with average-cost criterion, a team of agents sequentially interacts with the environment and observes local incentives. We focus on the case that the global reward is a sum of local rewards, the joint policy factorizes into agents' marginals, and full state observability. To date, few global optimality guarantees exist even for this simple setting, as most results yield convergence to stationarity for parameterized policies in large/possibly continuous spaces. To solidify the foundations of MARL, we build upon linear programming (LP) reformulations, for which stochastic primal-dual methods yields a model-free approach to achieve \emph{optimal sample complexity} in the centralized case. We develop multi-agent extensions, whereby agents solve their local saddle point problems and then perform local weighted averaging. We establish that the sample complexity to obtain near-globally optimal solutions matches tight dependencies on the cardinality of the state and action spaces, and exhibits classical scalings with respect to the network in accordance with multi-agent optimization. Experiments corroborate these results in practice.
\end{abstract}

%

\section{Introduction}\label{sec:intro}

In multi-agent reinforcement learning (MARL), a collection of agents repeatedly interact with their environment and are exposed to localized incentives. This framework has gained traction in recent years through successful application to autonomous vehicular networks \citep{wang2018reinforcement}, games \citep{vinyals2019grandmaster}, and various settings in econometrics \citep{tesauro2002pricing,lussange2021modelling}. At the core of MARL is a Markov Decision Process (MDP) \citep{puterman2014Markov}, which determines the interplay between agents, states, actions, and rewards. We focus on the standard objective, whereby the goal of the network agents is to discern policies so as to maximize the long-term accumulation of instantaneous rewards, which may be written as a node-separable sum of all localized rewards \citep{nedic2009distributed}.

Defining the team reward in this way implies that agents seek to cooperate towards a common goal, which may be contrasted with competitive or mixed settings \citep{bacsar1998dynamic}. Due to the surge of interest in MARL, disparate possible technical settings have been considered, which span how one defines MDP transition dynamics; the observability of agents trajectories, the availability of computational resources at a centralized location, and the protocol by which agents exchange information. We consider the case that agents have global knowledge of the state and action (in contrast to partial observability \citep{mahajan2016decentralized,krishnamurthy2016partially}, which necessitates pooling information as in centralized training decentralized execution (CTDE) \citep{foerster2016learning,leibo2017multi,foerster2017stabilising,rashid2018qmix}). Further, we hypothesize that the team's joint policy factorizes into the product of marginals, which is referred to as \emph{joint action learners} (JAL) \citep{claus1998dynamics,lee2020optimization}.

Our focus is on \emph{decentralized training} of JAL, which means agents' rewards and policy parameters are locally held and private. Numerous recent works on MARL operate in this setting, as in multi-agent extensions of temporal difference (TD) learning~\citep{lee2018stochastic,doan2019finite}, Q-learning~\citep{kar2013cal}, value iteration~\citep{wai2018multi,qu2019value}, and actor-critic~\citep{lowe2017multi,zhang2018fully}. In these works, agents may communicate according to the connectivity of a possibly time-varying graph, which is intimately connected to multi-agent optimization. \begin{table*}[htb!]
	\centering
	\resizebox{1.8\columnwidth}{!}	{\begin{tabular}{|c|c|c|c|}
			\hline
			References             & Rewards          &  Setting   & Sample Complexity 
			\\ 
			\hline
			\citep{raveh2019pac} &      Discounted              & Centralized  & $\widetilde{\mathcal{O}}\left(\frac{n|\mathcal{S}||\mathcal{A}|}{(1-\gamma)^2}\right)$ \\ 
			\hline
			\citep{wang2020randomized}, \citep{xu2020voting}&      Average              & Centralized, Parallel  & $ \Omega\left(\tau^2t_{mix}^2\frac{\mathcal{E}_{0} |\mathcal{S}||\mathcal{A}|}{\epsilon^2}\log\frac{1}{\delta}\right)$ \\ 
			\hline
			\citep{raveh2019pac} &      Average              & Parallel  & $\widetilde{\mathcal{O}}\left(\frac{n|\mathcal{S}||\mathcal{A}|}{(1-\gamma)^2}\right)$ \\ 
			\hline
			\citep{qu2020scalable2}  &      Average             &   Decentralized  &    ----
			\\ 
			\hline
		\textbf{	This work} &         Average            & Decentralized &       { $ \Omega\left(\tau^2t_{mix}^2\frac{\sqrt{n}\mathcal{E}_{0} |\mathcal{S}||\mathcal{A}|D(\Gamma, \rho)}{\epsilon^2}\log\frac{1}{\delta}\right)$}          \\ \hline
	\end{tabular}}
	\label{table}
	\caption{ MARL for average rewards case. The proposed scheme is the first decentralized algorithm in the average reward settings with PAC sample complexity. Here, $\epsilon$ is the accuracy parameter, $\delta$ denotes the high probability parameter, $n$ is the number of agents, $D(\Gamma, \rho):=\big[  \frac{1+\Gamma}{1-\rho}\big]$ where $\Gamma $$=$$ (1-\eta/4n^2)^{-2}$, $\rho $$=$$ (1-\eta/(4n^2))^{1/B}$, $B$ is the network strong connectivity parameter, and $\eta$ is a lower-bound on the entries of the mixing matrix.}
	\vspace{-0mm}
\end{table*}

 In the aforementioned references (limited to discounted objectives), convergence guarantees are mostly asymptotic, apply only to MARL sub-problems as policy evaluation (estimating the value function assuming a fixed policy \citep{sha2020fully,heredia2020finite}), or due to implied non-convexity induced by policy parameterization, cannot avoid spurious\footnote{Spurious here should be interpreted in the sense of stationary points that are far from global optimality.} policies \citep{qu2020scalable2,qu2020scalable} -- see \citep{zhang2020global} for further details. 

For these reasons, we focus on LP reformulation of RL in the average reward setting \citep{Ka83Thesis,Ka94I,de2003linear}, for which stochastic primal-dual method achieves \emph{optimal sample complexity} in the centralized tabular case \citep{wang2020randomized}. \footnote{In particular, we use $\mathcal{O}(1)$ to denote an absolute constant, and $\tilde{\mathcal{O}}(1)$ to hide polylog factors in $|\mathcal{S}|$, $|\mathcal{A}|$, and $\epsilon$, which are the respective cardinalities in the state and action spaces, and $\epsilon$ is some pre-defined optimization error.} Our goal is to understand which settings the
 $\tilde{\mathcal{O}}(\frac{\tau^2 t_{mix}^2 |\mathcal{S}||\mathcal{A}|}{\epsilon^2} \log (1/\delta))$
complexity achieved for finding an $\epsilon$-optimal solution with probability $1-\delta$ in the centralized case \citep{wang2020randomized}[Theorem 4] may be translatable to the multi-agent setting when agents may only exchange \emph{local} information with their one-hop neighbors. In particular, when agents combine localized stochastic primal-dual methods with weighted averaging to diffuse information across the network  \citep{nedic2009distributed,chen2012diffusion}, we seek to determine whether the optimal sample complexity of the LP approach generalizes to the average-reward tabular MARL. 
 Our contributions are to: 
 
  %
  %
  (i)  propose a novel multi-agent variant of the dual LP formulation of RL, where agents' decisions are defined by estimates of an average state-action occupancy measure and value vector, and consensus constraints are imposed on agents' localized estimates (Sec. \ref{sec:prob}). 
  
  %
  (ii) owing to node-separability of the Lagrangian relaxation of the resulting optimization problem, we derive a  decentralized model-free training mechanism based on a stochastic variant of primal-dual method that employs Kullback-Lieber (KL) divergence as its proximal term in the space of occupancy measures (Sec. \ref{sec:alg}), together with local weighted averaging.
  
  %
  (iii) establish that the number of samples required to attain near-globally optimal solutions matches tight dependencies on the cardinality of the state and action spaces \citep{wang2020randomized}, and exhibits classical scalings with the size of the team in prior theory \citep{nedic2009distributed}.
  
  %
  (iv) demonstrate the experimental merits of this approach in cooperative navigation problems.
  %

{\noindent \bf Additional Context.} Local averaging as a strategy for information mixing in multi-agent optimization is outperformed by schemes based upon Lagrange multiplier exchange, e.g., primal-dual method \citep{koppel2015saddle}, alternating direction method of multipliers (ADMM) \citep{boyd2011distributed}, and dual reformulations \citep{terelius2011decentralized}. In this work, however, we opt for a primal-only approach to enforcing consensus for simplicity and its compatibility with Perron-Frobenius theory \citep{chung1997spectral}. 

We further focus on the case where the communications network is a structural component of the problem setting, as in \citep{lowe2017multi,zhang2018fully}. However, a separate but related body of works estimate the communications architecture when agents' behavior is fixed using graph neural networks \citep{eccles2019biases,ahilan2020correcting,bachrach2020negotiating} or statistical tests for correlation between agents' local utilities~\citep{qu2020scalable,lin2020distributed}.

To the best of our knowledge, none of the aforementioned works deal with the average reward settings in MARL, with the exception of \citep{qu2020scalable2}. However, it provides asymptotic-only analysis. By contrast, the probably approximately correct (PAC) sample complexity results given here are unique to the MARL average-reward setting, and may be seen as a multi-agent generalization of \citep{wang2020randomized}. Critical to this generalization is a novel Lyapunov function that result to jointly tracks the convergence of the primal-dual iterates and the consensus error.  We note that PAC results have been developed for average-reward MARL in \citep{xu2020voting}. However, it operates under a setting where  policy and reward information are globally shared agents at each step, which in the optimization literature is known as \emph{parallel} \citep{bertsekas2015parallel}, not decentralized, as the updates cannot be executed with local and neighboring information only. For results most similar to this work in the MARL setting, please see Table \ref{table}.

\section{Problem Formulation}\label{sec:prob}
We consider MARL problems among agents who share a globally observable state, but take actions and observe rewards which are distinctly local. In this context, agents seek to coordinate in order to maximize the team's cumulative return of rewards, which is a sum over all locally observed rewards. More specifically, we consider a time-varying network $\mathcal{G}^t=(\mathcal{V},E^t, W^t)$ of $n$ agents $\mathcal{N}$$=$$\{1,2,\dots,n\}$, where agent $i\in\mathcal{N}$ may communicate with its neighbors $j$ if they share an edge $(i,j)\in{E}^t$, and no others, at a given time $t$. The weight matrix $W^t $$:=$$ [w_{ij}^t] \in \mathbb{R}^{n\times n}$, where $w_{ij}^t \geq 0$ and $w_{ij}^t = w_{ji}^t$ $\forall i,j,t$, assigns weights to each edge $(i,j)$. One canonical example of $w_{ij}^t$ is the relative degree between agent $i$ and $j$ at time $t$: $w_{ij}^t=d_i^t/(d_i^t+d_j^t)$, with $d_i^t$ as the degree, or number of nodes that are a one-hop neighbor of agent $i$.

With the network structure clarified, we now detail how the states, actions, and rewards interconnect. Precisely, at each time, agent $i \in \mathcal{V}$ observes the current system state $s \in \mathcal{S}$ and synchronously takes an action $a_i \in \mathcal{A}_i$, which is concatenated as the joint action $a $$=$$(a_1,...,a_n) \in \mathcal{A}_1\times\cdots\times\mathcal{A}_n$. The state space $\mathcal{S}$ and the constituent action spaces $\mathcal{A}_i$ are discrete finite sets with respectively $|\mathcal{S}$ and $|\mathcal{A}|$ elements. Trajectories are Markovian, that is, upon execution of the joint action $a$, the state transitions to next state $s'$ with probability $p_{s,s'}(a):=\mathbb{P}( s \mid s, a)$. That the joint action $a$ is observed by all agents after execution is needed to ensure full observability, i.e., that the MARL problem can be defined by an MDP. After the joint action $a$ is executed in state $s$, each agent $i$ receives a reward $r^i_{a,s}(s) \in [0,1]$, only known to the agent $i$. The system reward $r_{a,s}$ is defined as the aggregation of local rewards $r{(s,a)} $$:=$$\frac{1}{n} \sum_{i=1}^{n} r^i_{a,s}$.
The goal of the cooperative agents here is the maximization of the \emph{global} cumulative return defined as
\begin{align}\label{eq:main_prob}
    \max_{\pi} J_{\pi}(s) := \lim_{T \rightarrow \infty} \quad \frac{1}{T} \mathbb{E} \left[ \sum_{t=0}^{T-1} r_{a,s} \Big| s_0 = s \right],
\end{align}
where $\pi$ denotes the joint policy of all agents, that is, a probability distribution over joint action-space given system state, $\pi:\mathcal{S} \times \mathcal{A} \rightarrow [0,1]$. The joint policy prescribes the probability that a joint action $a=(a_1,\dots,a_n)$ is taken by the collection of the agents when in system state $s$, which we assume factors into marginals of each individual agent's policy: $\pi(a | s) $$=$$ \prod_{i=1}^{N} \pi_i(a_i | s)$. That is, the local policies are statistically independent, and are further denoted as  $\pi_i(a_i | s)$ which define the probability of taking action $a_i$ by agent $i$ when in state $s$. Moreover, the expectation in \eqref{eq:main_prob} is over the product measure associated with state transition dynamics and the policy known as the ergodic state occupancy measure.

Our specific goal in this work is the design of policy optimization schemes to solve \eqref{eq:main_prob} such that each agent, upon the basis of its local action selections and rewards, together with information exchange amongst neighbors, in possession of global state-action information, learns local policy parameters that result in the overall team attaining the optimal value \eqref{eq:main_prob}. We place specific emphasis upon the non-asymptotic convergence of such schemes and their scaling with respect to the parameters off network $\mathcal{G}$. Moreover, we consider in the model-free setting, i.e., the dynamics of the environment (the transition probabilities and transitional rewards) are unknown to the agents, but a simulation oracle is available to the team to generate state-action-reward tuples $(s,a,r)$.
We require that the transition dynamics for a fixed policy define an irreducible Markov chain: for each state pair $(s,s')$ and any policy $\pi$, there exist $t$ such that the probability that the system transitions from state $s$ to state $s'$ under policy $\pi$ in $t$ time-steps is non-zero. This condition is sufficient to ensure the limit in \eqref{eq:main_prob} exists and $J_{\pi}(s) =: \lambda_{\pi}$ for all states $s$. Equivalently, the average cost is independent of the initial state in the system. Further, the optimal policy is time-invariant.

 Towards transforming \eqref{eq:main_prob} into a workable form for deriving iterative model-free updates, we note that an optimal policy satisfies the \textit{average-cost Bellman equation} \citep{bertsekas1995dynamic} written as 
%
\begin{align} \label{eq:bellman_optimality}
   \!\! \lambda_{\pi}\!+\! v_s\! = \!\max_{a \in \mathcal{A}} \left \{ \sum_{s'} p_{s,s'}(a) r_{a,s} + \sum_{s'} p_{s,s'}(a) v_{s'} \right \},
\end{align}
for all $s \in \mathcal{S}$. Denote solutions to the Bellman's equation by pairs $(\lambda^*, v^*)$, where scalar $\lambda^*=\max J_\pi (s)$ in \eqref{eq:main_prob} is unique and equal to the optimal average cost. The value vector $v\in\mathbb{R}^s$ (which aggregates scalars $v_s$ for each $s\in\mathcal{S}$) is called a \textit{differential reward function} and is unique up to a constant. Uniqueness is imposed by $(v^*)^T\xi^* $$=$$ 0$, where $\xi^*$ is the stationary distribution under the optimal policy $\pi^*$, i.e. $P^{\pi^*} \xi^* $$=$$ \xi^*$. Note that each policy $\pi$ is associated with a transition probability $P_{s,s'}^{\pi} $$:=$$ \sum_{a} \pi(a|s) p_{s,s'}(a)$ for all $s,s' \in \mathcal{S}$, and a stationary state distribution $\xi^{\pi}$, which is a probability distribution that remains unchanged in the Markov chain as time progresses, i.e. $P^{\pi} \xi^{\pi}= \xi^{\pi}$. The differential reward function characterizes the transient effect of the initial state under a policy $\pi$. 

Continue then by noting that the optimal joint policy $\pi^*$ may be formulated as the following LP~\citep{DeFarias2003}:
\begin{equation} \label{dual_mat}
\begin{aligned}
   & \max_{\mu\in\mathbb{R}^{|\mathcal{S}|\times |\mathcal{A}|}}  \quad \sum_{a\in\mathcal{A}} \mu(a)^T r(a) \\
    &\textrm{s. t.\ } \begin{cases}
      \sum_{a\in\mathcal{A}} (I - P_a^T) \mu_a =0, \quad \forall s\\
      \sum_{s\in\mathcal{S},a\in\mathcal{A}} \mu(s,a) = 1
      \mu_{a,s} \geq 0 \quad \forall a,s
    \end{cases},
\end{aligned}
\end{equation}
where $I$ is an identity matrix of the appropriate size and $P_a \in \mathbb{R}^{|\mathcal{S}| \times |\mathcal{S}|}$ is the matrix whose $(s, s')$-th entry equals to $p_{s,s'}(a)$. For every feasible point of the above linear program $\mu = (\mu(a))_{a \in \mathcal{A}}$, the $\xi^{\pi} $$=$$ (\xi^{\pi}_s)_{s \in \mathcal{S}}$ is the stationary state distribution where $\xi^{\pi}_s $$=$$ \sum_{a} \mu_{a,s}$, and $\sum_{x,a} \mu(x,a) r_{a,s}$ corresponds to the average reward $\lambda_{\pi}$ of policy $\pi$ where $\pi(a|s) $$=$$ \frac{\mu_{a,s}}{\xi^{\pi}_s}$. 
Moreover, $\mu(a)\in\mathbb{R}^{|\mathcal{S}|}$ denotes the unnormalized occupancy measure over the state space $\mathcal{S}$ for each action $a\in\mathcal{A}$, whose stacking over the action space $\mathcal{A}$ is denoted as $\mu\in\mathbb{R}^{|\mathcal{S}|\times |\mathcal{A}|}$. Through normalization, one may recover the associated policy $\pi$ for any feasible $\mu$ as $    \xi^{\pi}_s $$=$$ \sum_{a\in\mathcal{A}} \mu(s,a)$, and $\pi(a|s)$$ =$$ \frac{\mu(s,a)}{\sum_{a\in\mathcal{A}} \mu(s,a)}$,
and from the definition of $\xi^{\pi}_s$ and $\pi(a|s)$, it holds that 
%
	$\mu(s,a) = \xi^{\pi}_s \pi(a|s)$. 
%
Then, an optimal joint policy $\pi^*$ can be constructed by normalizing the occupancy measures associated with the solution to the above linear program. See \citep{puterman2014Markov} and references therein for details. 
\begin{align} \label{opt_pi}
    \pi^*(a|s) = \frac{\mu^*(s,a)}{\sum_{a} \mu^*(s,a)}.
\end{align}
By substituting the definition of the global reward $r(s,a)$ in terms of the local rewards $r^i(s,a)$ into\eqref{dual_mat}, we obtain a {multi-agent optimization} problem with the global variables $\mu_{a,s}$ corresponding to joint policy $\pi$:
\begin{equation} \label{ma_lp}
\begin{aligned}
    \max_{\mu\in\mathbb{R}^{|\mathcal{S}|\times |\mathcal{A}|}} & \quad \sum_{i=1}^{n} \sum_{a\in\mathcal{A}} \mu(a)^T r^i(a) \\
    \textrm{subject to:}& \begin{cases}
      \sum_{a} (I - P_a^T) \mu(a) =0 \quad \forall s\in\mathcal{S} \\
      \sum_{s,a} \mu(s,a) = 1\\
      \mu(s,a) \geq 0 \quad \forall s\in\mathcal{S}, a\in\mathcal{A}
    \end{cases}.
\end{aligned}
\end{equation}
%
To solve \eqref{ma_lp}, agents must cooperate in their policy search. With each agent only exercising control over their localized policy, the globally optimal  joint policy $\pi^*(a|s)$ may be obtained via \eqref{opt_pi}. Specifically, under the previously mentioned independence assumption and knowledge of the state-action information, each agent may obtain its local policy by marginalizing the other agents' policies out of the the optimal joint policy as follows
%
%
 $   \pi_i(a_i|s) = \sum_{a_{-i}\in\mathcal{A}} \frac{\mu({s,(a_i,a_{-i})}}{\sum_{a}\mu({s,a})}$
where $a_{-i}$ denotes the joint action of all agents except the $i$-th agent; $a_i$ denotes the action associated with the $i$-th agent, and the joint action of all agents is denoted by $a$, i.e. $a$$=$$(a_i, a_{-i})$$=$$(a_1,\dots, a_n)$.
%
	\begin{algorithm2e}[htb!]
		\caption{Randomized Multi-agent Primal-dual (RMAPD) Algorithm}
		\label{alg_single}
		\textbf{Input:}$\epsilon > 0$, $\mathcal{S}$, $\mathcal{A}$, $t_{mix}^*$, $\tau$\\ 
		Set $v^i=0 \in \mathbb{R}^{|\mathcal{S}|}$, $\pi_i = \frac{1}{|\mathcal{A}_i|\mathbf{e}}\in \mathbb{R}^{|\mathcal{A}_i|}$, $\forall s \in \mathcal{S}$, $i \in \mathcal{N}$\\
		Set $T = {(\tau t_{mix}^*)}^2 |\mathcal{S}| |\mathcal{A}|$, $M=4t_{mix}^* + 1 $\\
		Set $\beta = \frac{1}{t^*_{mix}}\sqrt{\frac{\log(|\mathcal{S}\|\mathcal{A}|)}{2|\mathcal{S}\|\mathcal{A}|T}}$, $\alpha= |\mathcal{S}|t_{mix}^* \sqrt{\frac{\log(|\mathcal{S}\|\mathcal{A}|)}{2 |\mathcal{A}|T}}$	\\
		\For{iteration $t=0,1,2,...$}{
			\For{agent $i = 1,2,...,N$}{
				Observe the system state $s$, \\
				Execute action $a_i \sim \pi_i(\cdot|s)$ \\ 
				Observe local reward $r^i_{s,s'}(a)$ \\
				Send $\!(\mu^t_i,\!v_i^t)$ to $j\in n_i$, receive $\!(\mu^t_j,\!v_j^t)$.\\
				Compute local weighted averages [cf. \eqref{eq:consensus_round}]
				%
				$	\widetilde{\mu}_i^t = \sum_{j=1}^{n} w_{ij}^t \mu_j^t$, 
				$	\widetilde{v}_i^t = \sum_{j=1}^{n} w_{ij}^t v_j^t \;$,
			\\
				%
				%
				%
				%
				%
				%
				%
				%
				%
				\! Conduct entropic ascent w.r.t. $\!\!\mu_i^{t+1}\!$  in \eqref{eq:mu_updates}:
							\begin{align*}
								\!\!\!\!\!\mu_i^{t + \frac{1}{2}}(s,a) &= \frac{\widetilde{\mu}_i^t(s,a) \exp( \Delta_i^{t+1}(s,a))}{\sum_{s'} \sum_{a'} \widetilde{\mu}_i^t(s,a) \exp( \Delta_i^{t+1}(s',a'))}\\
								\mu_i^{t+1} &= \argmin_{\mu_i \in \mathcal{U}} D_{KL} (\mu_i \| \mu_i^{t + \frac{1}{2}}), 
							\end{align*}
				with dual gradient $\Delta_i^{t+1}(s,a)$ in \eqref{eq:dual_gradient}.\\
				Update value vector for agent $i$ via \eqref{eq:value_updates} as
				%
								%
							$	v_i^{t+1} = \Pi_{\mathcal{V}} [\widetilde{v}_i^t +  d_i^{t+1}]$
				with $d_i^{t+1}$ in \eqref{eq:primal_gradient}
			}	
		}
	\end{algorithm2e}

With the setting clarified, we next shift to developing a decentralized model-free algorithm to solve \eqref{eq:main_prob} upon the basis of Lagrangian relaxation.

\section{Randomized Primal-Dual Method}\label{sec:alg}
In this section, we reformulate the multi-agent LP of \eqref{ma_lp} as a saddle point problem by considering its Lagrangian relaxation. 
%
%
In particular, we formulate the following saddle point problem
\begin{align}\label{eq:saddle_point_problem}
\!\!\!\!\!\!\   \min_{v \in \mathcal{V}} \max_{\mu \in \mathcal{U}}  L(\mu, v)\! :=\! \sum_{i=1}^{n} \sum_{a \in \mathcal{A}}\! \mu(a)^T[(P_a-I)v\!+\!r^i(a)].\!\!
\end{align}
Note that we have computed the transpose of the constraint to simplify the expression. Under Assumptions \ref{tau-stationary} and \ref{mixing_time} introduced in Sec. \ref{sec:convergence}, we may establish that the primal-dual optimal pair $(v^*,\mu^*)$ of \eqref{eq:saddle_point_problem} belong to the following restricted setsfor the value $\mathcal{V}\subset \mathbb{R}^{|\mathcal{S}|}$ and occupancy measures $\mathcal{U}\subset \mathbb{R}^{|\mathcal{S}|\times|\mathcal{A}|}$ defined as
\begin{align}\label{eq:feasible_sets}
\mathcal{V} =&  \{v \in \mathbb{R}^{|\mathcal{S}|} \Big| \quad \|v\|_{\infty}\ \leq 2t_{mix}  \} \;, \quad
%
\\
\mathcal{U} =& \big \{\mu \!=\! (\mu_a)_{a \in \mathcal{A}} \Big| \mathbf{e}^T\mu \!=\! 1,  \mu\!\geq\!0, \sum_{a \in \mathcal{A}} \mu(a) \geq \frac{1}{\sqrt{\tau}|\mathcal{S}|}\mathbf{e} \big \},\nonumber
\end{align}
where $t_{mix}$ is the mixing time of the Markov chain which characterizes how fast the Markov decision process reaches its stationary distribution from any state under any policy (Assumption \ref{mixing_time}), and $\tau$ is a constant greater than one which characterizes how much the stationary distribution varies as the policy varies (Assumption \ref{tau-stationary}). The definitions of these feasible sets is borne out of the analysis, and mirrors \citep{wang2020randomized}.
Next, we note that the Lagrangian of \eqref{eq:saddle_point_problem} is node-separable. Specifically, by defining the local Lagrangian for agent $i \in V$ as 
\begin{equation}\label{eq:local_lagrangian}
L_i(\mu, v) := \sum_{a}\mu(a)^T [(I - P_a)  v + r^i(a)],
\end{equation}
where $v$ is a column vector with $v_s$ as its $s$-th component, then the Lagrangian of the multi-agent problem $    L(\mu, v) $ may be decomposed into a sum over local Lagrangian $L_i(\mu, v)$ as
%
    $L(\mu, v) = \sum_{i=1}^{n} L_i(\mu, v)$,
%
which permits us to simplify the saddle point problem as
\begin{align}\label{eq:local_lagrangian2}
    \min_{v \in \mathcal{V}} \max_{\mu \in \mathcal{U}} \quad L(\mu, v) &= \sum_{i=1}^{n} L_i(\mu, v).
\end{align}
The min-max problem in \eqref{eq:local_lagrangian2} is convex in $v$ and concave in $\mu$. We note that the variables $v$ and $\mu$ are common among all the agents in the network and we are interested in solving the problem in a distributed manner. This expression in \eqref{eq:local_lagrangian2} is suggestive of employing a solution methodology upon the basis of a decentralized stochastic primal-dual method, which is the focus of the following subsection.
\subsection{Stochastic Primal-Dual Method}\label{sec:alg_iterates}
We propose applying stochastic primal-dual method to solve \eqref{ma_lp}, which, owing to the node-separability of the Lagrangian, yields a decentralized scheme for policy optimization. 
%
%
In particular, in order to solve the saddle point problem, we note that agents must access estimate the global reward, but they lack access. Instead, agents only observe local rewards. To address this issue, we allow each agent to track a distinctly localized estimate $v_i\in\mathcal{S}$ of the value, which are substituted in place of the global value vector $v$ in \eqref{eq:saddle_point_problem}, and similarly with respect to the occupancy measure $\mu_i^t$, which are in lieu of the global primal-dual pair $(\mu^t, v^t)$.
Then, agent $i$ cooperates with other agents through a weighted averaging of its primal and the dual variables, i.e. a convex combination $\widetilde{\mu}_i^t$ (resp. $\widetilde{v}_i^t$) of its own estimate $\mu_i^t$ (resp. $v_i^t$) with the estimates received from those of its neighbors $j\in n_i$ at time $t$:
\begin{align}\label{eq:consensus_round}
		\widetilde{\mu}_i^t = \sum_{j=1}^{n} w_{ij}^t \mu_j^t, \qquad 
		\widetilde{v}_i^t = \sum_{j=1}^{n} w_{ij}^t v_j^t \;,
\end{align}
Then, each agent takes a gradient descent (respectively, ascent) step to minimize (respectively, maximize) the local Lagrangian function $L_i$,  followed by a projection onto the constraint set $\mathcal{U}$ (respectively, $\mathcal{V}$). However, since the transition dynamics model is unavailable to agent $i$ (in the form of transition matrix $P_a$), it cannot to evaluate the constraint in \eqref{dual_mat}. This precludes the evaluation of primal and dual gradients of the Lagrangian, which necessitates stochastic approximations of these quantities, which we present jointly with respective step-size parameters $\beta$ and $\alpha$ as
\begin{align}\label{eq:dual_gradient}
    &\hat{\nabla}_{\mu_i}L_i = \Delta_i^{t+1} =  \beta \frac{{ v_i^t(s')- v_i^t(s)}+r_i^t(s,s',a)-M}{\widetilde{\mu}_i^t(s,a)}.\mathbf{e}_{s,a}, 
    \nonumber
    \\
    &\qquad \qquad\qquad \qquad\textrm{with probability $\widetilde{\mu}_i^t(s,a)$}
\end{align}
\begin{align}\label{eq:primal_gradient}
    &\hat{\nabla}_{v_i}L_i = d_i^{t+1}  ={\frac{{\mu}_i^t(s,a)}{\widetilde{\mu}_i^t(s,a)}}\alpha (e_s-e_{s'}), 
    \\
    &\qquad \qquad\qquad \qquad \textrm{with probability $\widetilde{\mu}_i^t(s,a)p_{s,s'}(a)$}, \nonumber
\end{align}
where $M := 4t_{mix}+1$ is a ``shift parameter" which ensures sufficient decrease of a certain martingale process defined in terms of the KL divergence that arises in the analysis (to be made precise later), and the superscript $t$ denotes the value of the variable at time $t$. Moreover $t_{mix}$ is the mixing time of the Markov chain induced by a fixed policy (Assumption \ref{mixing_time}). Here $\mathbf{e}_{s,a}$ is the indicator variable which is $1$ for $(s,a)$ and null otherwise. Further, $e_s$ denotes the standard basis vector with $1$ in slot $s$ and null otherwise. Note that we adopt the convention that, at time-step $t$, the variables with superscript $t$ are known, and the superscript $t+1$ indicates an update direction in terms of random variables realized at time $t$. An additional point of note is that the gradient with respect to the value vector is $-\alpha (e_{s'}-e_{s})$, which we swap to cancel out the negative. Then, using these update directions, stochastic primal-dual method is such that at every $t \geq 0$, each agent $i$ generates new estimates $\mu_i^{t+1}$, $v^i_{t+1}$ as
\begin{align}\label{eq:mu_updates}
    &\mu_i^{t+1} = \argmin_{\mu_i \in \mathcal{U}} D_{KL} (\mu_i \| \mu_i^{t + \frac{1}{2}}), 
    \\
    & \ \textrm{where} \quad \mu_i^{t + \frac{1}{2}}(s,a) = \frac{\widetilde{\mu}_i^t(s,a) \exp( \Delta_i^{t+1}(s,a))}{\sum_{s'} \sum_{a'} \widetilde{\mu}_i^t(s,a) \exp( \Delta_i^{t+1}(s',a'))}\nonumber
    \\
     &   v_i^{t+1} = \Pi_{\mathcal{V}} [\widetilde{v}_i^t +  d_i^{t+1}]\label{eq:value_updates},
\end{align}
where $\Pi_{\mathcal{V}}$ is a Euclidean projection onto the set $\mathcal{V}$, and $d_i^{t+1}$ is given in \eqref{eq:primal_gradient}. Moreover, $\Delta^{t+1}_t$ is the gradient of the local Lagrangian with respect to $\mu_i$ in \eqref{eq:dual_gradient}. Note that the update on $\mu$ is mirror-ascent with a Kullback-Leibler (KL) divergence over the unnormalized probability simplex centered at $\widetilde{\mu}_i^t$, whereas the gradient step on the value vector $v$ is a simple projected gradient descent centered at $\widetilde{v}_i^t$. The descent step on $v$ is written in terms of an addition due to the cancellation of a negative, as mentioned after \eqref{eq:primal_gradient}. We assume algorithm initialization as $\mu_i = 0$ and $v_i = 0$ for all $i \in \mathcal{N}$. The overall MARL policy optimization scheme based upon randomized primal-dual solutions to the LP formulation is summarized as Algorithm \ref{alg_single}.
%


\section{Convergence Analysis}\label{sec:convergence}
In this section, we establish the non-asymptotic convergence of the proposed algorithm in the sense that agents' local primal-dual variables (a) achieve consensus and (b) converge to the primal-dual optimal pair of their local Lagrangians \eqref{eq:local_lagrangian}. As a consequence, upon the basis of local observations and information exchange with neighbors, agents are able to solve \eqref{ma_lp}, and hence \eqref{eq:main_prob}. We divide the analysis of the algorithm in two steps. First, we establish that all local estimates achieve consensus. Second, we show that the consensus vectors are in fact a pair of primal-dual optimal solution. To establish these results, we state some conditions are required on the graph $\mathcal{G}^t$ next.
%
\begin{assumption} \label{strong_connectivity}[Strong Connectivity]
There exists a positive integer $B$ such that graph $(\mathcal{N},
\cup_{l=0}^{B-1} E_{t+l})$ is strongly-connected for any $t \geq 0$, i.e., every node is reachable from another in at most $B$ time-steps.
\end{assumption}
%
\begin{assumption} \label{mixing_matrix}
For all $i \in V$ and $t \geq 0$: 
(a) there exists a scalar $\eta \in (0 < \eta < 1)$ such that $w_{ij}^t \geq \eta$ when $j \in \mathcal{N}_i^t$ , and $w_{ij}^t = 0$ otherwise\label{mixing_lower_bound}; 
(b) $\sum_{j=1}^{n} w^t_{ij} = \sum_{i=1}^{n} w^t_{ij} $$=$$ 1$; that is, the mixing matrix $W^t$ is doubly stochastic \label{mixing_stochastic}.
\end{assumption}
Assumption \ref{strong_connectivity} ensures that after a union of $B$ time-slots, the network is connected, which ensures information propagates across the network. Assumption \ref{mixing_matrix} ensures that an agent sufficiently balances the weighting of its own information with that of other agents. Assumption \ref{mixing_matrix} ensures that the mixing matrices have a  Perron-Frobenius eigenvalue associated with an eigenvector whose entries are all $1$, i.e., the existence of a vector satisfying consensus.
We also make two assumptions on the MDP stationary distribution and mixing time of the chain:
\begin{assumption} \label{tau-stationary}[Ergodic Decision Process]
The Markov decision process is $\tau$-stationary in the sense that it is ergodic under any stationary policy $\pi$ and there exists $\tau>1$ such that
$    \frac{1}{\sqrt{\tau}|\mathcal{S}|}\mathbf{e} \leq \xi^{\pi} \leq \frac{\sqrt{\tau}}{|\mathcal{S}|}\mathbf{e}$,
where $\mathbf{e}$ is a vector of all $1$'s.
\end{assumption}
\begin{assumption} \label{mixing_time} [Fast-Mixing Markov Chains]
The Markov decision process is $t_{mix}$-mixing in the sense that
 $   t_{mix} \geq \max_{\pi} \min \big \{ t\geq 1~\Big|~ \| (P^{\pi})^t(s,.) - \xi^{\pi}\|_{TV} \leq \frac{1}{4}, \forall s \in \mathcal{S} \big \}$,
where $\|.\|_{TV}$ is the total variation norm.
\end{assumption}
The factor $\tau$ characterizes the variability of the stationary distribution with respect to the policy. $t_{mix}$ defines how fast the MDP reaches its stationary distribution from any state under any policy $\pi$.
\subsection{Primal-dual Optimality}\label{subsec:optimality}
To show that the consensus vector coincides with a pair of primal-dual optimal solution, we show that, at each iteration of the algorithm, the local iterates get closer to the local primal-dual optimal pair in expectation. It turns out that to do so, we must first show that agents' estimates reach approximate consensus, and then construct a Lyapunov function with respect to quantities defined in terms of network averages. We proceed to doing so next,

{\bf \noindent Achieving consensus.} We establish that the local iterates converge to the global mean at a specified rate in terms of the lower bound on the mixing weights, the diameter of the network, and the strong connectivity parameter. The consensus error must be characterized for both the value vector and occupancy measure estimates, which motivate the following network-aggregated averages at time $t$:
\begin{align}\label{eq:average_vector}
	\overline{\mu}^{t} = \frac{1}{n}\sum_{i=1}^{n} \mu_i^{t}, \  \quad
	\overline{v}^{t} =  \frac{1}{n} \sum_{i=1}^{n} v_i^{t} \; .
\end{align}
It also turns out to be convenient to define the auxiliary sequence
\begin{align}\label{eq:auxiliary_sequence}
	%
	q_i^{t} := \Pi_{\mathcal{V}} [\widetilde{v}_i^t +  d_i^{t+1}] - \widetilde{v}_i^t ,
\end{align}
where $q_i^{t}$ represent the error between the weight-averaged iterates $\widetilde{\mu}_i^t$ (resp. $\widetilde{v}_i^t$) and their previous update following  projection/composition with a proximal operator. We note that a similar technique for minimization problems is considered in \citep{chen2021distributed}, but here we are considering a different minimax setting, necessitating analyzing the consensus error in both the primal and dual variables. See Lemma \ref{dist_mean} in the appendix for the analysis  of consensus error. 

{\bf \noindent Lyapunov Function Construction.} Next we define a decrement process that tracks the evolution of the averaged primal and dual iterates to the primal-dual optimal pair, which eventuates in our ability to formalize the overall convergence rate of Algorithm \ref{alg_single}. Consider the Lyapunov function $\mathcal{E}_{t}$ and duality gap quantifier $\mathcal{D}_t$, defined as
\begin{align}\label{eq:lypunov_function}
	\mathcal{E}_{t}&:=\frac{1}{n}\sum_{i=1}^n D_{KL}(\mu^* \| \mu_i^{t})  +\frac{1}{2|\mathcal{S}|t_{mix}^2}\left\| \overline{v}^{t}- v^*\right\|^2 \nonumber \\
	\mathcal{D}_t&:=\lambda^* + \sum_{a \in A}\left[\overline{\mu}^t(a)^T[( I-P_a) v^* +r_a] \right].
\end{align}
The first term of $\mathcal{E}_{t}$ quantifies the sum of $KL$ divergences between the optimal $\mu^*$ and local occupancy measures $\mu_i^{t}$, and the second term quantifies the sub-optimality of the average sequence $\overline{v}^{t}$. In $\mathcal{D}_t$, we track the constraint violation of \eqref{ma_lp}.
\begin{lemma}\label{lemma:KL_divergence_martingale}
	With Lyapunov function $\mathcal{E}_{t}$ and duality gap quantifier $\mathcal{D}_t$ in\eqref{eq:lypunov_function}, the iterates of Algorithm \ref{alg_single} exhibits approximate stochastic descent:
	\begin{align}\label{eq:KL_divergence_martingale}
		&\mathbb{E}\left[\mathcal{E}_{t+1}\mid \mathcal{F}_t \right]
		\leq \mathcal{E}_{t} -\beta\mathcal{D}_t +\beta^2\widetilde{\mathcal{O}}\left(|\mathcal{S}||\mathcal{A}|t_{mix}^2\right) 	
\\
		&\hspace{0cm}+\frac{\beta}{n}\!\sum_{i=1}^{n}\!\sum_{a \in A}[ \left(\overline{v}^t -v_i^t \right)^T\left( (I - P_a)^T(\widetilde{\mu}_i^t(a)- \mu^*(a))\right)]  
\nonumber
\\
&\hspace{0.2cm}+ \frac{\beta}{n}\sum_{i=1}^{n}\sum_{a \in A}[\left(\widetilde{\mu}_i^t(a)-\overline{\mu}^t(a)\right)^T[(P_a - I) v^* +r_a] ],\nonumber
	\end{align}
	where $\beta$ denotes the constant step size. 
\end{lemma}
See Appendix \ref{proof_lemma:KL_divergence_martingale} for proof. This result is a generalization of \cite[Proposition 9]{wang2020randomized} to multi-agent settings. For $n=1$, the last two terms are null, which simplifies to the proposition in the aforementioned reference. In generalizing it to the multi-agent setting, we note that existing analyses of multi-agent stochastic optimization methods based on consensus protocol rely on finite variance conditions \citep{li2018stochastic}. However, the dual gradient to be evaluated at consensus variable $\widetilde{\mu}_i^t(a)$ is required for unbiasedness in the gradient evaluation in  \eqref{eq:dual_gradient}, may cause unbounded noises to the stochastic gradient estimates. An additional complication is the joint treatment of consensus error in primal and dual variables owing to the structure of the minimax objective \eqref{eq:saddle_point_problem}, which is not treated in any of the earlier works \citep{boyd2011distributed,nedic2009distributed,kar2013cal,li2018stochastic,chen2021distributed,sha2020fully}.

Next, we present the main result of this subsection which upper bounds the duality gap as  follows. 
\begin{theorem}\label{eq:KL_divergence_martingale_final400}
	For the time-averaged sequence of occupancy measures $\hat{{\mu}}$ $=$ $\frac{1}{T}\sum_{t=0}^{T-1}\frac{1}{n}\sum_{i=1}^{n} \mu_i^{t}$, after $T$ number of iterations of Algorithm \ref{alg_single}, with the step size selection $\beta=\mathcal{\widetilde{\mathcal{O}}}\left(\sqrt{\frac{\mathcal{E}_{0}}{{ \sqrt{n} |\mathcal{S}||\mathcal{A}| \tilde{t}^2_{mix}D(\Gamma, \rho)}T}}\right)$, it holds that 
	\begin{align}\label{eq:KL_divergence_martingale_final4200}
			\mathbb{E}&\left[\sum_{a \in A}\left[[ v^*-P_av^*  +r_a]^T\hat{{\mu}}(a) \right]\right] +\lambda^*	
		\nonumber
		\\	&\hspace{1.2cm}	\leq \widetilde{\mathcal{O}}\left( \tilde{t}_{mix}\sqrt{\frac{\sqrt{n}\mathcal{E}_{0} |\mathcal{S}||\mathcal{A}|D(\Gamma, \rho)}{T}}\right),
	\end{align}
where $\tilde{t}_{mix}=1+4{t}_{mix}$, $D(\Gamma, \rho):=\big[  \frac{1+\Gamma}{1-\rho}\big]$ such that $\Gamma $$=$$ (1-\eta/4n^2)^{-2}$ and $\rho $$=$$ (1-\eta/(4n^2))^{1/B}$, $B$ is the network strong connectivity parameter.
\end{theorem}
\begin{figure*}[htb!]
	\centering
	\subfigure[ \scriptsize Grid world.]{
		\includegraphics[scale = 0.28]{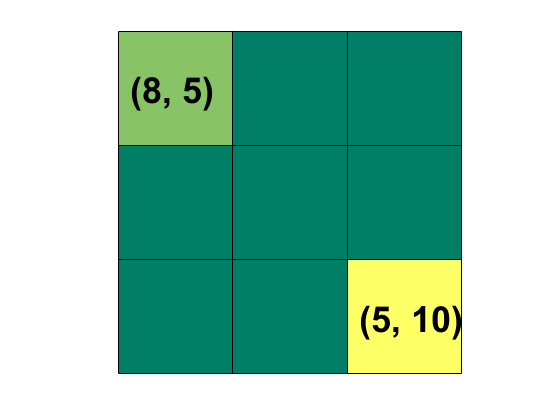}
		\label{fig:regret_rs_6}
	}
	\subfigure[  \scriptsize $M$$=$$3$, $n$$=$$2.$]{
			\includegraphics[scale = 0.15]{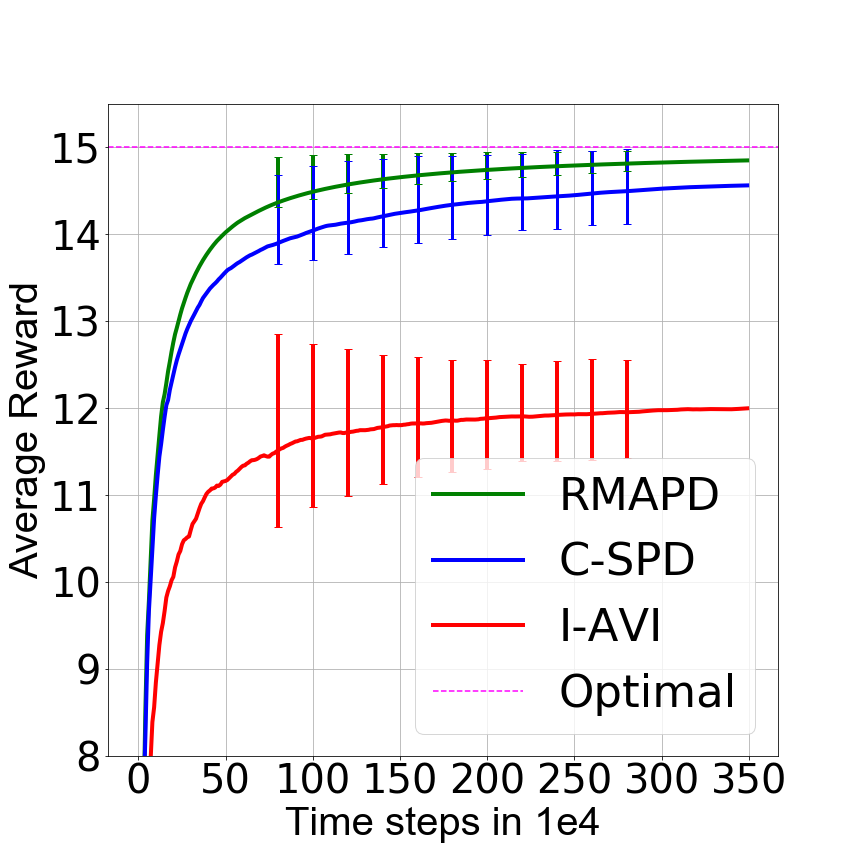}	
		\label{fig:regret_rs_12}
	}
\subfigure[ \scriptsize  $M$$=$$2$, $n$$=$$3$. ]{
\includegraphics[scale = 0.15]{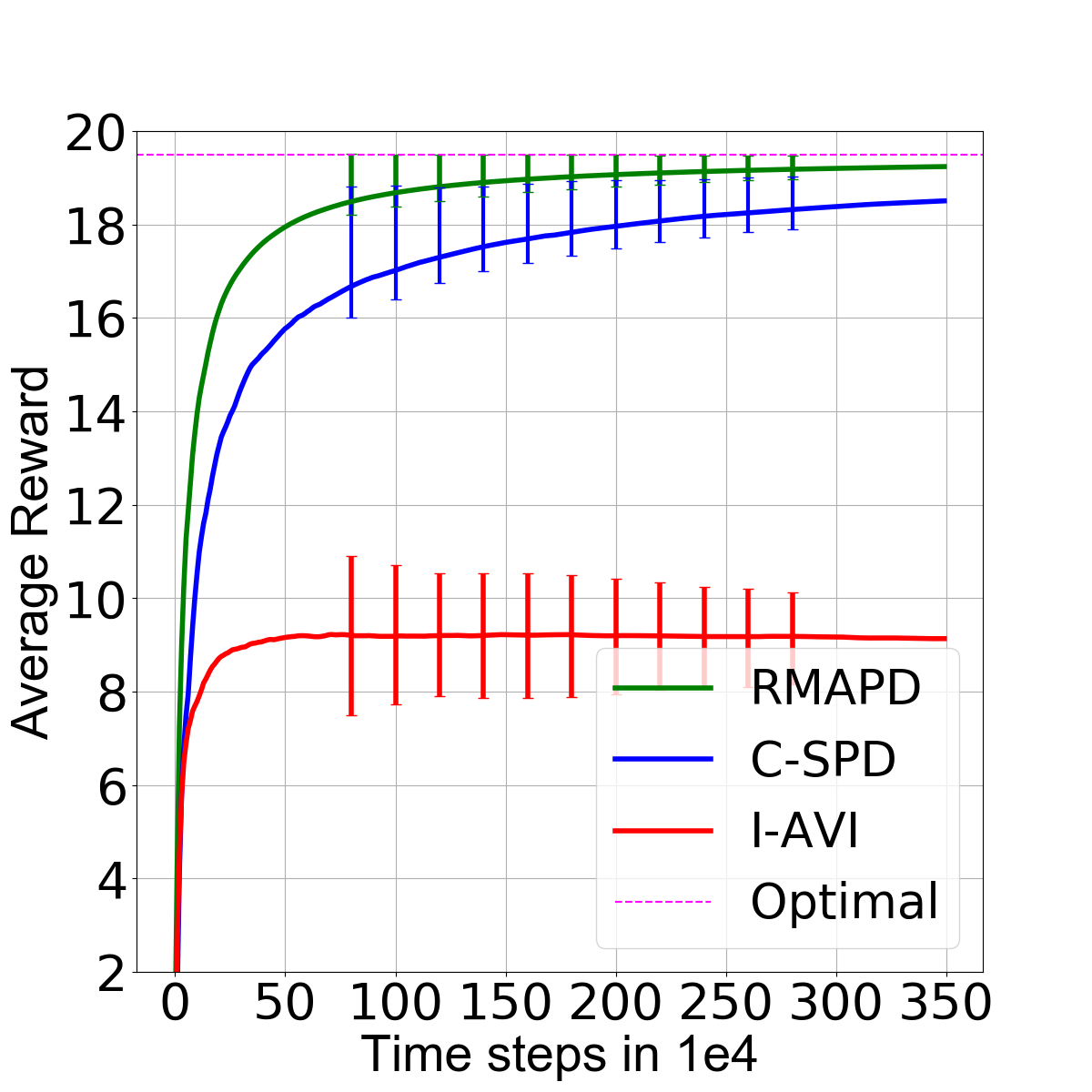}
		\label{fig:regret_gw}
	}    
	\caption{ We compare Algorithm \ref{alg_single} with its centralized counterpart (centralized stochastic primal-dual, or C-SPD), and with independent approximate value iteration (I-AVI). Fig. \ref{fig:regret_rs_6} shows the grid world environment for the experiments. Fig. \ref{fig:regret_rs_12} compares the average reward for all the mentioned algorithms. It shows that RMAPD is able to learn the optimal policy equivalent to centralized technique. We note that I-AVI fails to learn the optimal policy because agents are not cooperating with each other. Fig. \ref{fig:regret_gw} shows similar result for $M=2, n=3$. }
	\label{fig:per_agent_regret}
	\vspace{0in}
\end{figure*}
See Appendix \ref{proof_eq:KL_divergence_martingale_final400} for proof. Observe that the duality gap characterization is nonstandard from typical saddle point problems \citep{nedic2009subgradient}, that is, the left-hand side of \eqref{eq:KL_divergence_martingale_final4200} characterizes how the averaged dual variable $\hat{{\mu}}$ evaluated at the constraint at the optimal primal variable $(v^*,\lambda^*)$, and not an additional presence of the dual sub-optimality. This is a special structural consequence of the LP setting that breaks down for general nonlinear objectives or constraints. This upper bound characterizes the number of times the complementary slackness condition is violated on average, which facilitates deriving the sample complexity required to achieve an $\epsilon$-optimal policy with high probability. Next we shift focus to this result.

\subsection{From Duality Gap to Average Reward}
We first derive the convergence in probability result for the proposed algorithm in next Lemma \ref{prop}. 
 \begin{lemma}\label{prop}
 	Suppose Algorithm \ref{alg_single} is run for $T=\Omega\left(\tau^2\tilde{t}_{mix}^2\frac{\sqrt{n}\mathcal{E}_{0} |\mathcal{S}||\mathcal{A}|D(\Gamma, \rho)}{\epsilon^2}\right)$ iterations. Then it outputs a policy  $\hat\pi$=$\frac{1}{T}\sum_{t=1}^{T}\overline\pi^t$ such that $\lambda^*-\epsilon 	\leq\lambda_{\hat\pi}$, 
 	with probability $2/3$, meaning, we output an $\epsilon$ optimal policy with probability $2/3$.
 	 \end{lemma}
See Appendix \ref{prob_1} for proof. Lemma \ref{prop} establishes that Algorithm \ref{alg_single} converges to $\epsilon$-optimal policy with probability $2/3$.  To boost the success probability to near $1$, we develop a strategy where one runs Algorithm \ref{alg_single} multiple times and selects the best outcome. This procedure is formalized in Algorithm \ref{alg_single2}. With this meta-strategy in practice, we may establish that one can indeed achieve achieve comparable sample complexity to Lemma \ref{prop} but with high probability.

  \begin{algorithm2e}[h]
  	\caption{Meta-Randomized Multi-agent Primal-dual (M-RMAPD) Algorithm}
  	\label{alg_single2}
  	\textbf{Input:}$\epsilon > 0$, $\mathcal{S}$, $\mathcal{A}$, $t_{mix}^*$, $\tau$\\ 
   	Run the Algorithm \ref{alg_single} for $K$ number of iterations with precision $\frac{\epsilon}{3}$ and denote the output as $\overline{\pi}^{(1)}, \cdots, \overline{\pi}^{(K)}$. \\
 For each output policy $\overline{\pi}^{(k)}$, conduct the approximate value evaluation for $L=\tilde{\mathcal{O}}\left(\frac{t_{mix}}{\epsilon^2}\log\left(\frac{4K}{\delta}\right)\right)$ time steps and obtain $\overline{Y}^{(k)}$ which is approximate value evaluation with precision level $\epsilon/3$ and prob. $\frac{\delta}{2K}$. \\
Output $\widetilde\pi=\overline\pi^{(k^*)}$ such that $k^*= \argmax_{k}\overline{Y}^{(k)}$. 
  \end{algorithm2e}
  \begin{theorem}\label{thm:main_sample_complexity}
  Under Assumption \ref{strong_connectivity}-\ref{mixing_time}, if we run Algorithm \ref{alg_single2} for $K=\log_{1/3} \left(\frac{\delta}{2}\right)$,   then we output an approximate policy $\widetilde\pi$ such that $\lambda_{\widetilde\pi} \geq  \lambda^*-{\epsilon}$ with at least probability $1-\delta$. Hence, the total number of samples required are given by $T=\Omega\left(\tau^2\tilde{t}_{mix}^2\frac{\sqrt{n}\mathcal{E}_{0} |\mathcal{S}||\mathcal{A}|D(\Gamma, \rho)}{\epsilon^2}\cdot\log\frac{1}{\delta}\right)$ We define $D(\Gamma, \rho):=\big[  \frac{1+\Gamma}{1-\rho}\big]$ such that $\Gamma $$=$$ (1-\eta/4n^2)^{-2}$ and $\rho $$=$$ (1-\eta/(4n^2))^{1/B}$ where $B$ is the network strong connectivity parameter.
  	\end{theorem}

  See Appendix \ref{theorem_2} for proof. To the best of our knowledge, the result in Theorem \ref{alg_single} is the first to characterize the sample complexity of MARL schemes with high probability to achieve global optimality. We accentuate that we are able to discern explicit dependence upon the mixing time $t_{mix}$ and network parameters with tight dependence upon the cardinalities of the state and action spaces. 
 %


\section{Experiments} \label{sec_experiments}
In this section, we evaluate the practical merit of the proposed algorithm for MARL. In particular, since we are interested in cooperative multi-agent RL, we consider an experimental  setting where the need for aggregating the policy learnt by individual agents via consensus or centralized training is important. 
Hence, we consider a cooperative navigation problem in a $M\times M$ grid world environment shown in Fig. \ref{fig:regret_rs_6}. Each agent is equipped with action $\mathcal{A}_i=\{\uparrow, \rightarrow, \downarrow,\leftarrow\}$ and observe $M\times M$ grid as the local state space $\mathcal{S}_i$.  In this environment, each agent receives a reward $r_i$, and the common goal is to reach a state with maximum average reward across the agents. For instance, 
%
%
in our grid world environment as depicted in \ref{fig:regret_rs_6} for $M=3$ and $n = 2$, the agent $1$ and agent 2 receive a reward of $8$ and $5$ in the top left grid, and a reward of $5$ and $10$ in the lower rightmost grid, respectively, when they reach there simultaneously, and zero otherwise. This settings ensures that the cooperative behavior would result in higher average reward as compared to a non-cooperative behavior. 


We solve the grid world navigation problem using the proposed RMAPD algorithm and present the average cumulative reward returns in Fig. \ref{fig:per_agent_regret}.  We compare the performance of the proposed decentralized algorithm with a centralized LP solver and also a variant of approximate value iteration where each agent operates independently of all others to maximize its local average reward. The plot in Fig. \ref{fig:regret_rs_12} shows that the proposed algorithm iterates converge to the centralized optimal solution and is significantly better than the independent learning schemes for a grid of size $M=3$ with $n=2$ agents, and similarly for $M=2,n=3$ in Fig. \ref{fig:regret_gw}.
%
%
%

In the experiments, we run $20$ independent iterations of all the algorithms for $10^7$ timesteps and plot the average rewards. We observe that convergence is reached sooner in RMAPD as compared to centralized training because of the need to explore a larger state space per agent before converging as opposed to its decentralized counterpart. Further, the importance of the consensus mechanism of the proposed algorithm is highlighted by the much lower average reward achieved by the independent approximate value iteration.

\section{Conclusions}
In this work, we considered a multi-agent reinforcement learning problem with average-reward criterion. The problem has been studied in literature in single agent settings only to date, for which randomized LP solvers achieve optimal PAC bounds in terms of the cardinality of the state and action spaces. We generalized such approaches to multi-agent settings by combining randomized LP solvers with consensus averaging, and elucidated their PAC bounds for the same. Interestingly, by the use of a novel Lyapunov function in the convergence analysis, the dependence upon the state and action space cardinality is still maintained similar to centralized counterpart, with an additional dependence on the way information propagates across the multi-agent network.

As a future direction, we will develop variants of this framework that can operate with parameterized occupancy measures and differential value vectors, such that the scaling is with respect to the parameterization rather than the state and action spaces.

\bibliographystyle{icml2020}
\bibliography{main,bibliography}
\onecolumn

\appendix
\newpage

\section*{\centering {Supplementary Material for \\ ``Convergence Rates of  Average-Reward Multi-agent \\Reinforcement Learning via Randomized Linear Programming"}}

\section{Proof of Lemma \ref{lemma:KL_divergence_martingale}} \label{proof_lemma:KL_divergence_martingale}

Before starting the proof of Lemma \ref{lemma:KL_divergence_martingale}, we provide the intermediate lemmas which are useful in the analysis as follows. 

We proceed to defining the transition matrix $\Phi(t,s)$, for all $(t,s)$ with $t \geq s \geq 0$ as a product of weight matrices, given by
\begin{align} \label{transition_matrix}
	\Phi(t,s) = W^t W^{t-1}\dots W^s,
\end{align}
as the product of weight matrices associated with the time-varying graph $\mathcal{G}^t=(\mathcal{V},E^t,W^t)$.
The following is a key fact regarding the transition matrices \eqref{transition_matrix}. 
\begin{proposition} \label{tran_bound}\citep{nedic2009distributed}[Proposition 1]
	Under Assumptions \ref{strong_connectivity} and \ref{mixing_matrix}, for all $i, j \in \mathcal{V}$ and all $(t, s)$ with $t \geq s \geq 0$, we have
	\begin{align}
		\Big|[\Phi(t,s)]_{ij}-\frac{1}{n}\Big| \leq \Gamma \rho^{t-s}
	\end{align}
	where $\Gamma $$=$$ (1-\eta/4n^2)^{-2}$ and $\rho $$=$$ (1-\eta/(4n^2))^{1/B}$, with, $n$ the number of nodes in $\mathcal{N}$, $B$ is the strong-connectivity parameter of Assumption \ref{strong_connectivity}, and $\eta$ the lower-bound on the weights $w_{ij}$ in Assumption \ref{mixing_matrix}.
\end{proposition}
Observe that for a static graph, we have $\Phi(t,s) = W^{t-s+1}$ and $B=1$, that is, the weighting matrix is constant over time and the graph is assumed to be strongly connected.
Next, we establish that the iterates have bounded deviation from the mean.
\begin{lemma} \label{dist_mean}
	Let Assumptions \ref{strong_connectivity} and \ref{mixing_matrix} hold. Recall constants $\rho$ and $\Gamma$ from Lemma \ref{tran_bound}.
	
	(i) Denote as $\{v_i^{t}\}_{t \geq 0}$ the sequence of differential value vector estimates generated by Algorithm \ref{alg_single}. Then, for constant step size $\alpha$, for all $i \in V$ and $t \geq 0$, we have
	\begin{align*}
		\mathbb{E}[\|\mathbf{v}^{t} \!-\! (\frac{1}{n}\mathbf{e}\mathbf{e}^T\otimes I_{|\mathcal{S}|}) \mathbf{v}^t\|\!\mid\!\mathcal{F}_t] \!\leq \!\mathcal{O}(\sqrt{n}\alpha) \cdot{Z_t}(\Gamma,\rho),
	\end{align*}
	where ${Z_t}(\Gamma,\rho)$$:=$$[1$$+$$  \frac{\Gamma(1-\rho^{t-1})}{1-\rho}]$,~$\mathbf{v}^t$$=$$[[v_1^t]^T ;  \cdots ; [v_n^t]^T]\in\mathbb{R}^{n|\mathcal{S}|}$ stacks value vector estimates $v_i^t$, and $\left(\frac{1}{n}\mathbf{e}\mathbf{e}^T\otimes I_{|\mathcal{S}|}\right) \mathbf{v}^t$ stacks of $n$ copies of average vector $\overline{v}^t$ [cf. \eqref{eq:average_vector}] of  dimension $n|\mathcal{S}|$, and $ I_{|\mathcal{S}|}\in\mathbb{R}^{|\mathcal{S}|\times |\mathcal{S}|}$ is the identity matrix.

	(ii) The sequence of occupancy measure estimates $\{\mu_i^{t}\}_{t \geq 0}$ generated by Algorithm \ref{alg_single} under constant step size $\alpha$, satisfy for all $i \in V$ and $t \geq 0$
	\begin{align}\label{eq:dist_mean2}
		\mathbb{E}&\left[\|	\boldsymbol{\mu}^{t} - \left(\frac{1}{n}\mathbf{e}\mathbf{e}^T\otimes I_{|\mathcal{S}||\mathcal{A}|}\right) 	\boldsymbol{\mu}^{t} \| \mid \mathcal{F}_t\right] 
		\\
		&\leq (2(4t_{mix}+1)\beta\sqrt{|\mathcal{S}| |\mathcal{A}|} \sqrt{n})\cdot{Z_t}(\Gamma,\rho),\nonumber
	\end{align}
	where ${Z_t}(\Gamma,\rho)$$:=$$[1$$+$$  \frac{\Gamma(1-\rho^{t-1})}{1-\rho}]$,~$	\boldsymbol{\mu}^{t}$$=$$[[\mu_1^t]^T ;  \cdots ; [\mu_n^t]^T]\in\mathbb{R}^{n|\mathcal{S}||\mathcal{A}|}$ stacks occupancy measure estimates $\mu_i^t$ across all $i$, and $\left(\frac{1}{n}\mathbf{e}\mathbf{e}^T\otimes I_{|\mathcal{S}||\mathcal{A}|}\right) \boldsymbol{\mu}^{t}$ stacks of $n$ copies of average vector $\overline{\mu}^t$ [cf. \eqref{eq:average_vector}] of  dimension $n|\mathcal{S}||\mathcal{A}|$, and $ I_{|\mathcal{S}||\mathcal{A}|}\in\mathbb{R}^{|\mathcal{S}||\mathcal{A}|\times |\mathcal{S}||\mathcal{A}|}$ is the identity matrix.
\end{lemma}
See Appendix \ref{proof_dist_mean} for proof. This result is employed in the analysis of the evolution of the localized primal-dual iterations \eqref{eq:mu_updates} - \eqref{eq:value_updates} when we decompose the sub-optimality with respect to the global saddle point problem \eqref{eq:saddle_point_problem} into a consensus error and an optimization error especially to address the consensus error.

Next, we prove the intermediate results which leads to the convergence of Lemma \ref{lemma:KL_divergence_martingale}.  
%
In particular, we start with the KL Divergence difference $ D_{KL}(\mu^* \| \mu_i^{t+1})-D_{KL}(\mu^* \| \mu_i^{t}) $ in the following Lemma \ref{lemma:KL_divergence_contraction}.

\begin{lemma}\label{lemma:KL_divergence_contraction}
	Consider the KL Divergence-based average between the local dual variable iterates defined in \eqref{eq:mu_updates} and the optimal occupancy measure $\mu^*$ defined by \eqref{eq:local_lagrangian} satisfies the following approximate decrement:
	\begin{align}\label{eq:KL_divergence_contraction}
		\frac{1}{n}\sum_{i=1}^nD_{KL}(\mu^* \| \mu_i^{t+1})-\frac{1}{n}\sum_{i=1}^nD_{KL}(\mu^* \| \mu_i^{t})  	\leq &\frac{1}{n}\sum_{i=1}^{n} \sum_{s\in\mathcal{S}}\sum_{a\in\mathcal{A}}\left(\widetilde{\mu}_i^t(s,a)-{\mu}^*(s,a)\right) \Delta_i^{t+1}(s,a) 
		\nonumber
		\\ &+ \frac{1}{2n}\sum_{i=1}^{n} \sum_{s\in\mathcal{S}}\sum_{a\in\mathcal{A}}  \widetilde{\mu}_i^t(s,a) (\Delta_i^{t+1}(s,a))^2. 
	\end{align}
\end{lemma}
The proof of Lemma \ref{lemma:KL_divergence_contraction} is provided in Appendix \ref{proof_lemma:KL_divergence_contraction}. {In comparison tot the analysis for single agent settings in \citep{wang2020randomized}, the right hand side of \eqref{eq:KL_divergence_contraction} now depends upon the consensus variable $\widetilde{\mu}_i^t$.} Next, we establish the conditional expectation of update direction of the local occupancy $\Delta_i^{t+1}$ with respect to a difference of the averaged occupancy measure with respect to the optimal admits an expression in terms of the constraint violation of \eqref{ma_lp}.
\begin{lemma}\label{lemma:occupancy_decrement}
	For arbitrary $s \in S$ and $a \in A$, we have
	\begin{align}\label{eq:occupancy_decrement}
		\frac{1}{n}\sum_{i=1}^{n}\sum_{s \in S} \sum_{a \in A} (\widetilde{\mu}_i^t(s,a) - \mu^*(s,a)) \mathbb{E}[\Delta_i^{t+1}(s,a)\mid \mathcal{F}_t] 
		&=  \frac{\beta}{n}\sum_{i=1}^{n}\sum_{a \in A} (\overline{\mu}^t(a) - \mu^*(a))^T[(P_a - I) v_i^t +r_a].
	\end{align}
\end{lemma}
The proof of Lemma \ref{lemma:occupancy_decrement} is provided in Appendix \ref{proof_lemma:occupancy_decrement}.
\begin{lemma}\label{lemma:mean_square_dual}
	The $\widetilde{\mu}_i^t(s,a) $ probability-weighted conditional mean-square of the local occupancy measure update $\Delta_i^{t+1}(s,a)$ satisfies the following boundedness condition:
	\begin{align}\label{eq:mean_square}
		\sum_{s \in S} \sum_{a \in A} \widetilde{\mu}_i^t(s,a) \mathbb{E}[(\Delta_i^{t+1}(s,a))^2\mid\mathcal{F}_t] &\leq  4(4t_{mix}+1)^2\beta^2|\mathcal{S}| |\mathcal{A}| \; .
		%
	\end{align}    
\end{lemma}
The proof of Lemma \ref{lemma:mean_square_dual} is provided in Appendix \ref{proof_lemma_mean_Swuare_dual}. 
\begin{lemma}\label{lemma:value_decrement}
	Consider the norm difference between the average differential value vector update defined by \eqref{eq:value_updates} and the global-optimal $v^*$ as the minimizer of $\sum_i L_i(\mu, v)$ in \eqref{eq:local_lagrangian}. This difference satisfies the following recursion:
	\begin{align}\label{eq:value_decrement}
		\mathbb{E}\left[\left\| \overline{v}^{t+1}- v^*\right\|^2 \mid \mathcal{F}_t\right]& \leq
		\left\| \overline{v}^t - v^* \right \|^2 + \frac{2 \alpha}{n} \left(\overline{v}^t - v^* \right)^T\left( \frac{1}{n}\sum_{i=1}^n \sum_{a\in\mathcal{A}} (I - P_a)^T\widetilde{\mu}_i^t(a)\right)  + \frac{1}{n} \mathcal{O}(\alpha^2)
	\end{align}
\end{lemma}
The proof of Lemma \ref{lemma:value_decrement} is provided in Appendix \ref{proof_lemma:value_decrement}. The statement of Lemma \ref{lemma:value_decrement} establishes the norm distance of average dual variable $\overline{v}^{t+1}$ from the optimal dual variable $v^*$.

The proof of Lemma \ref{lemma:KL_divergence_contraction}-\ref{lemma:value_decrement} is used to proof the statement of Lemma \ref{lemma:KL_divergence_martingale}. Begin by considering the expectation of \eqref{eq:KL_divergence_contraction} conditional on filtration $\mathcal{F}_t$. Then, employ Lemma \ref{lemma:occupancy_decrement} and Lemma \ref{lemma:mean_square_dual} for the second two terms on the right-hand side:
\begin{align}\label{eq:KL_divergence_martingale_start}
	\frac{1}{n}\sum_{i=1}^n \mathbb{E}\left[D_{KL}(\mu^* \| \mu_i^{t+1}) \mid \mathcal{F}_t \right] 	
	\leq & \frac{1}{n}\sum_{i=1}^nD_{KL}(\mu^* \| \mu_i^{t}) 
	+\frac{1}{n}\sum_{i=1}^{n} \sum_{s\in\mathcal{S}}\sum_{a\in\mathcal{A}}\left(\widetilde{\mu}_i^t(s,a)-{\mu}^*(s,a)\right)  \mathbb{E}\left[\Delta_i^{t+1}(s,a)\mid \mathcal{F}_t \right] 
	\nonumber
	\\
	&+ \frac{1}{2n}\sum_{i=1}^{n} \sum_{s\in\mathcal{S}}\sum_{a\in\mathcal{A}}  \widetilde{\mu}_i^t(s,a)\mathbb{E}\left[(\Delta_i^{t+1}(s,a))^2 \mid \mathcal{F}_t\right]
	\nonumber
	\nonumber
	\\
	\leq &\frac{1}{n}\sum_{i=1}^nD_{KL}(\mu^* \| \mu_i^{t}) + {}\frac{\beta}{n}\sum_{i=1}^{n}\sum_{a \in A} (\widetilde{\mu}_i^t(a) - \mu^*(a))^T[(P_a - I) v_i^t +r_a)] 
	\nonumber
	\\
	&+ 4(4t_{mix}+1)^2\beta^2|\mathcal{S}| |\mathcal{A}|.
\end{align}
%
Now, consider the expression for the value function decrease in Lemma \ref{lemma:value_decrement}. 
Multiply both sides of this inequality by $\frac{1}{2|\mathcal{S}|t_{mix}^2}$ (where $|\mathcal{S}|$ denotes the cardinality of the state space, and $t_{mix}$ is the mixing time of the MDP defined in Assumption \ref{mixing_time}) and add to \eqref{eq:value_decrement_mixing2} to obtain
\begin{align}\label{eq:KL_divergence_martingale_start2}
	\mathbb{E}\left[\frac{1}{n}\sum_{i=1}^n D_{KL}(\mu^* \| \mu_i^{t+1}) +\frac{1}{2|\mathcal{S}|t_{mix}^2}\left\| \overline{v}^{t+1}- v^*\right\|^2\mid \mathcal{F}_t \right]
	\leq &\frac{1}{n}\sum_{i=1}^n D_{KL}(\mu^* \| \mu_i^{t}) +\frac{1}{2|\mathcal{S}|t_{mix}^2}\left\| \overline{v}^{t}- v^*\right\|^2
	\nonumber
	\\
	&\quad  + \frac{\beta}{n}\sum_{i=1}^{n}\sum_{a \in A} (\widetilde{\mu}_i^t(a) - \mu^*(a))^T[(P_a - I) v_i^t +r_a)]  	\nonumber
	\\
	&\quad\quad  + \frac{2 \alpha}{n}\frac{1}{2|\mathcal{S}|t_{mix}^2}\!\! \left(\overline{v}^t\! -\! v^* \right)^T\!\!\left(\!\! \frac{1}{n}\sum_{i=1}^n \sum_{a\in\mathcal{A}} (I\! -\! P_a)^T\widetilde{\mu}_i^t(a)\!\!\right)	\nonumber
	\\
	&\quad\quad\quad +\frac{1}{2|\mathcal{S}|t_{mix}^2} \mathcal{O}(\alpha^2)+ 4(4t_{mix}+1)^2\beta^2|\mathcal{S}| |\mathcal{A}|.
\end{align}
Notice that the last two terms are respectively of order $\alpha^2$ and $\beta^2$, with a contrast in the order of the dependence on problem-dependent constants $|\mathcal{S}|$, $|\mathcal{A}|$, and $t_{mix}$. These terms may be judiciously balanced via selecting  ${\alpha= |\mathcal{S}|t^2_{mix}\beta}$, which yields
\begin{align}\label{eq:KL_divergence_martingale_start22}
	\mathbb{E}\left[\frac{1}{n}\sum_{i=1}^n D_{KL}(\mu^* \| \mu_i^{t+1}) +\frac{1}{2|\mathcal{S}|t_{mix}^2}\left\| \overline{v}^{t+1}- v^*\right\|^2\mid \mathcal{F}_t \right]
	\leq &\frac{1}{n}\sum_{i=1}^n D_{KL}(\mu^* \| \mu_i^{t}) +\frac{1}{2|\mathcal{S}|t_{mix}^2}\left\| \overline{v}^{t}- v^*\right\|^2
	\nonumber
	\\
	&\quad  + \frac{\beta}{n}\sum_{i=1}^{n}\sum_{a \in A} (\widetilde{\mu}_i^t(a) - \mu^*(a))^T[(P_a - I) v_i^t +r_a)]  	\nonumber
	\\
	&\quad\quad  + \frac{\beta}{n}\left(\overline{v}^t - v^* \right)^T\left( \frac{1}{n}\sum_{i=1}^n \sum_{a\in\mathcal{A}} (I - P_a)^T\widetilde{\mu}_i^t(a)\right)	\nonumber
	\\
	&\quad\quad\quad +\beta^2\widetilde{\mathcal{O}}\left(|\mathcal{S}||\mathcal{A}|t_{mix}^2\right).
\end{align}
where as defined earlier, $\widetilde{\mathcal{O}}$ denotes order-dependence that ignores polylog factors.
After grouping the third and fourth terms on the right-hand side of the previous expression, we can write
\begin{align}\label{eq:KL_divergence_martingale_start33}
	\mathbb{E}\left[\mathcal{E}_{t+1}\mid \mathcal{F}_t \right]
	\leq &\mathcal{E}_{t} +\frac{\beta}{n}\sum_{i=1}^{n}\sum_{a \in A}\left[(\widetilde{\mu}_i^t(a) - \mu^*(a))[(P_a - I) v_i^t +r_a] + \left(\overline{v}^t - v^* \right)^T\left( (I - P_a)^T\widetilde{\mu}_i^t(a)\right)\right]  	\nonumber
	\\
	& +\beta^2\widetilde{\mathcal{O}}\left(|\mathcal{S}||\mathcal{A}|t_{mix}^2\right).
\end{align}
where we have used the definition $\mathcal{E}_{t}:=\frac{1}{n}\sum_{i=1}^n D_{KL}(\mu^* \| \mu_i^{t+1})  +\frac{1}{2|\mathcal{S}|t_{mix}^2}\left\| \overline{v}^{t}- v^*\right\|^2$. Now, let us study second term on the right-hand side of the previous expression \eqref{eq:KL_divergence_martingale_start33} as
%
\begin{align}\label{eq:complimentary_slackness_analysis}
	\frac{\beta}{n}&\sum_{i=1}^{n}\sum_{a \in A}\left[(\widetilde{\mu}_i^t(a) - \mu^*(a))^T[(P_a - I) v_i^t +r_a] + \left(\overline{v}^t - v^* \right)^T\left( (I - P_a)^T\widetilde{\mu}_i^t(a)\right)\right] \nonumber
	\\
	& =	{\frac{\beta}{n}}\sum_{i=1}^{n}\sum_{a \in A}\left[(\widetilde{\mu}_i^t(a) - \mu^*(a))^T[(P_a - I) v_i^t +r_a] + \left(\overline{v}^t - v^* \right)^T\left( (I - P_a)^T(\widetilde{\mu}_i^t(a)- \mu^*(a))\right)\right]  
\end{align}
where we have used the dual feasibility of $\mu^*$, which means that $\sum_{a} (I - P_a^T) \mu^*(a) =0$. Now, add and subtract $v_i^t$ inside the parenthesis associating the value function difference between $\overline{v}^t$ and $ v^*$ in the preceding expression to obtain
\begin{align}\label{eq:complimentary_slackness_analysis2}
	\frac{\beta}{n}\sum_{i=1}^{n}&\sum_{a \in A}\left[(\widetilde{\mu}_i^t(a) - \mu^*(a))^T[(P_a - I) v_i^t +r_a] + \left(v_i^t- v^*+\overline{v}^t -v_i^t \right)^T\left( (I - P_a)^T(\widetilde{\mu}_i^t(a)- \mu^*(a))\right)\right]
	\nonumber\\
	& =\frac{\beta}{n}\sum_{i=1}^{n}\sum_{a \in A}\left[(\widetilde{\mu}_i^t(a) - \mu^*(a))^T[(P_a - I) v_i^t +r_a] + \left(v_i^t- v^*\right)^T\left( (I - P_a)^T(\widetilde{\mu}_i^t(a)- \mu^*(a))\right)\right]
	\nonumber\\
	&\qquad +	\frac{\beta}{n}\sum_{i=1}^{n}\sum_{a \in A}\left[ \left(\overline{v}^t -v_i^t \right)^T\left( (I - P_a)^T(\widetilde{\mu}_i^t(a)- \mu^*(a))\right)\right]
	\nonumber\\
	& =	\frac{\beta}{n}\sum_{i=1}^{n}\sum_{a \in A}\left[(\widetilde{\mu}_i^t(a) - \mu^*(a))^T[(P_a - I) v^* +r_a] \right]  \nonumber\\
	&\qquad +	\frac{\beta}{n}\sum_{i=1}^{n}\sum_{a \in A}\left[ \left(\overline{v}^t -v_i^t \right)^T\left( (I - P_a)^T(\widetilde{\mu}_i^t(a)- \mu^*(a))\right)\right] 
\end{align}
where the first equality uses the definition of the transpose, and the later regroups terms. We proceed to analyze both terms on the right-hand side of the preceding expression step by step. The former one is analogous to \citep{Wang2020}[Lemma A.6], but the later term is a novel instantiation of the consensus error due to decentralized computations. Unsurprisingly, Lemmas \ref{tran_bound} is useful to address it.

First, we focus on the first term on the right-hand side of \eqref{eq:complimentary_slackness_analysis2}. To do so, we exploit the linear complementarity of $(v^*,\mu^*)$ in \eqref{eq:saddle_point_problem}, i.e.,
$$\mu^*(s,a)[(P_a - I)v^* + r_a - \lambda^*\mathbf{e}]_s=0\; \text{ for all } s\in\mathcal{S}\; , $$
where $\lambda^*$ is defined as the optimal objective following \eqref{eq:bellman_optimality}, $\mathbf{e}$ is the vector in $\mathbb{R}^{|\mathcal{S}|}$ whose entries are all $1$. This fact can be substituted in place of $\mu^*(s,a)(P_a - I)v^* + r_a$ in first term on the right-hand side of \eqref{eq:complimentary_slackness_analysis2} as

\begin{align}\label{eq:complimentary_slackness_analysis3}
	\frac{\beta}{n}\sum_{i=1}^{n}\sum_{a \in A}\left[(\widetilde{\mu}_i^t(a) - \mu^*(a))^T[(P_a - I) v^* +r_a] \right] &=\frac{\beta}{n}\sum_{i=1}^{n}\sum_{a \in A}\left[\widetilde{\mu}_i^t(a)^T[(P_a - I) v^* +r_a] \right] - \sum_{a\in\mathcal{A}}\lambda^* \mu^*(a)\mathbf{e}\nonumber \\
	&=\frac{\beta}{n}\sum_{i=1}^{n}\sum_{a \in A}\left[\left(\overline{\mu}^t(a)-\overline{\mu}^t(a)+\widetilde{\mu}_i^t(a)\right)^T[(P_a - I) v^* +r_a] \right] - \lambda^* \nonumber \\
	&=- {\beta}\underbrace{(\lambda^* + \sum_{a \in A}\left[\overline{\mu}^t(a)^T[( I-P_a) v^* +r_a] \right])}_{:=\mathcal{D}_t}
	\nonumber
	\\
	&\qquad- \frac{\beta}{n}\sum_{i=1}^{n}\sum_{a \in A}\left[\left(\overline{\mu}^t(a)-\widetilde{\mu}_i^t(a)\right)^T[(P_a - I) v^* +r_a] \right]
\end{align}
Now, let us substitute the definition of $\mathcal{D}_t$ in \eqref{eq:complimentary_slackness_analysis3} into the right-hand side of \eqref{eq:KL_divergence_martingale_start33} together with the decomposition of the complimentary slackness in \eqref{eq:complimentary_slackness_analysis2}:
\begin{align}\label{eq:KL_divergence_martingale_final}
	\mathbb{E}\left[\mathcal{E}_{t+1}\mid \mathcal{F}_t \right]
	\leq& \mathcal{E}_{t} -\beta\mathcal{D}_t +\beta^2\widetilde{\mathcal{O}}\left(n|\mathcal{S}||\mathcal{A}|t_{mix}^2\right) + \underbrace{\frac{\beta}{n}\sum_{i=1}^{n}\sum_{a \in A}\left[ \left(\overline{v}^t -v_i^t \right)^T\left( (I - P_a)^T(\widetilde{\mu}_i^t(a)- \mu^*(a))\right)\right]}_{I}  
	\nonumber
	\\
	&+ \underbrace{\frac{\beta}{n}\sum_{i=1}^{n}\sum_{a \in A}\left[\left(\widetilde{\mu}_i^t(a)-\overline{\mu}^t(a)\right)^T[(P_a - I) v^* +r_a] \right]}_{II}.
\end{align}
which is as stated in Lemma \ref{lemma:KL_divergence_martingale}. $\hfill \square$.

\section{Proof of Lemma \ref{dist_mean}}\label{proof_dist_mean}

\subsection{Proof of Lemma \ref{dist_mean} Statement (i)} Begin by employing the definition of the the auxiliary sequences $\{p_i^t\}_{t\geq 0}$ and $\{q_i^t\}_{t\geq 0}$ for all agents $i \in \mathcal{V}$. First, we analyze $q_i^t$ and show its magnitude is upper-bounded. In particular, note that
\begin{align}\label{eq:upper_bound_q}
	\| q_i^t\| &=\| \Pi_{\mathcal{V}} [\widetilde{v}_i^t +  d_i^{t+1}] - \widetilde{v}_i^t\| \nonumber \\
	&\leq  \| \widetilde{v}_i^t +  d_i^{t+1} - \widetilde{v}_i^t\|
\end{align}
by the non-expansive property of the projection. Then, we may cancel out $\widetilde{v}_i^t$ and employ Lemma \ref{lemma:mean_square_dual} regarding the mean-square magnitude of $d_i^{t+1}$ to conclude:
\begin{align}\label{eq:upper_bound_mean_q2}
	\mathbb{E}[\| q_i^t\| \mid \mathcal{F}_t] \leq \mathbb{E}[\| d_i^{t+1}\| \mid \mathcal{F}_t] {\leq \mathcal{O}(\alpha)}.
\end{align}
Now, first we establish the bound on $\mathbb{E}[\|d_i^{t+1}\|^2\mid\mathcal{F}_t]$ as follows. For the analysis of the magnitude of $d_i^{t+1}$, we obtain by direct evaluation of its mean conditional on the filtration $\mathcal{F}_t$:
{$$\mathbb{E}[d_i^{t+1}\mid\mathcal{F}_t] = \alpha \sum_{a\in\mathcal{A}} (I - P_a)^T{\mu}_i^t(a)$$}
from the definition of the probability of the update direction in \eqref{eq:primal_gradient} and the local Lagrangian in \eqref{eq:local_lagrangian}. Now, observe that $\widetilde{\mu}_i^t$ may not satisfy the ergodicity constraint of Assumption \ref{tau-stationary} because weighted averaging [cf. \eqref{eq:consensus_round}] does not preserve feasibility, i.e., it may be outside dual feasible set $\mathcal{U}$ in \eqref{eq:feasible_sets}. However, since the weighting matrices are doubly stochastic (Assumption \ref{mixing_matrix}), $\widetilde{\mu}_i^t$ still belongs to the probability simplex, meaning its sum is unit. Together with the nonnegativity of the probability, means that the the inner product on the right-hand side of the preceding expression sums up to less than one, i.e., 
\begin{equation}\label{eq:value_update_bound}
	{ \alpha \sum_{a\in\mathcal{A}} (I - P_a)^T{\mu}_i^t(a)\leq \mathcal{O}(\alpha)}
\end{equation}
Next, in order to obtain the bound on $\mathbb{E}[\|d_i^{t+1}\|^2\mid\mathcal{F}_t]$, let us consider the expression as 
\begin{align}
	\mathbb{E}[\|d_i^{t+1}\|^2\mid\mathcal{F}_t]=& \mathbb{E}\left[\Bigg\|{\frac{{\mu}_i^t(s,a)}{\widetilde{\mu}_i^t(s,a)}}\alpha (e_s-e_{s'})\Bigg\|^2\mid\mathcal{F}_t\right]
	\nonumber
	\\
	\leq& \alpha ^2\mathbb{E}\left[\|e_s-e_{s'}\|^2\mid\mathcal{F}_t\right]
	\nonumber
	\\
	\leq& 4\alpha ^2,
\end{align}
{which establishes that $\mathbb{E}[\|d_i^{t+1}\|^2\mid\mathcal{F}_t]\leq \mathcal{O}(\alpha^2)$.} From the Jensen's inequality, we known that $f(\mathbb{E}(x))\leq \mathbb{E} f((x))$ for a convex function $f(x)$. Now, by selecting convex function $f(x)=x^2$, we can write $\left(\mathbb{E}[\|d_i^{t+1}\|\mid\mathcal{F}_t]\right)^2\leq \mathbb{E}[\|d_i^{t+1}\|^2\mid\mathcal{F}_t] \leq \mathcal{O}(\alpha^2)$, which further implies that $\mathbb{E}[\|d_i^{t+1}\|\mid\mathcal{F}_t]\leq \mathcal{O}(\alpha)$. Hence, we get \eqref{eq:upper_bound_mean_q2}.  
which we may stack over all $i\in\mathcal{V}$ to write $\mathbb{E}[\| q_t\| \mid \mathcal{F}_t] \leq \mathcal{O}(\sqrt{n}\alpha)$, where $\mathbf{q}_t=[[q_1^t]^T ;  \cdots ; [q_n^t]^T]\in\mathbb{R}^{n|\mathcal{S}|}$ stacks $q_i^t\in\mathbb{R}^{|\mathcal{S}|}$ [cf. \eqref{eq:auxiliary_sequence}] over all $i$. For the moment, consider that the value vector is scalar for each agent. Then, using this definition of $q_i^t$, we may write the the value vector $\mathbf{v}^t=[[v_1^t]^T ;  \cdots ; [v_n^t]^T]\in\mathbb{R}^{n|\mathcal{S}|}$ stacking $v_i^t$ across all $i$ as
\begin{align}\label{eq:stacked_value_mixing}
	\mathbf{v}^{t+1} &= \widetilde{\mathbf{v}}^t + \mathbf{q}^t\nonumber \\
	&= (W^t \otimes I_{|\mathcal{S}|}) \mathbf{v}^t + \mathbf{q}^t
\end{align}
where  $\otimes\mathcal{D}_t $ denotes the Kronecker product and $I_{|\mathcal{S}|}\in\mathbb{R}^{|\mathcal{S}|\times |\mathcal{S}|}$ is the $|\mathcal{S}|\times |\mathcal{S}|$ identity matrix. 
With these observations, shift to defining the respective sequences that track the distance to the mean [cf. \eqref{eq:average_vector}]
\begin{equation}\label{eq:consensus_error_vector}
	\delta_i^t = v_i^t - \overline{v}^t \; . 
\end{equation}
We first analyze the evolution of $\delta_i^t$ stacked over all $i$ defined as $\mathbf{\delta}^t=[[\delta_1^t]^T ;  \cdots ; [\delta_n^t]^T]\in\mathbb{R}^{n|\mathcal{S}|}$. Observe that $\delta^t$ may be defined in terms of the uniform matrix of weights all set to $(1/n)$, i.e.,
\begin{equation}\label{eq:consensus_error_vector_identity}
	\mathbf{\delta}^t = \mathbf{v}^t - \left(\frac{1}{n}\mathbf{e}\mathbf{e}^T\otimes I_{|\mathcal{S}|}\right) \mathbf{v}^t = \left[I_{n|\mathcal{S}|} - \left(\frac{1}{n}\mathbf{e}\mathbf{e}^T\otimes I_{|\mathcal{S}|}\right)  \right] \mathbf{v}^t
\end{equation}
where $\mathbf{e}$ is the vector of all $1$'s in $\mathbb{R}^{n}$, meaning that $\frac{1}{n}\mathbf{e}\mathbf{e}^T$ is an $n\times n$ matrix whose entries are all equal to $1/n$. Further, $I_{|\mathcal{S}|}\in\mathbb{R}^{|\mathcal{S}|\times |\mathcal{S}|}$ is the identity matrix, and similarly for $I_{n|\mathcal{S}|}$ in dimension $n|\mathcal{S}| \times n|\mathcal{S}| $ . Observe that $\delta^t$ satisfies the recursion:
\begin{align}\label{eq:consensus_error_recursion}
	\delta^{t+1} &=\left[I_{n|\mathcal{S}|} - \left(\frac{1}{n}\mathbf{e}\mathbf{e}^T\otimes I_{|\mathcal{S}|}\right)  \right] \mathbf{v}^{t+1}\nonumber \\
	&=\left[I_{n|\mathcal{S}|} - \left(\frac{1}{n}\mathbf{e}\mathbf{e}^T\otimes I_{|\mathcal{S}|}\right) \right] \left[(W^t \otimes I_{|\mathcal{S}|}) \mathbf{v}^t + \mathbf{q}^t\right]\nonumber \\
	&= (W^t \otimes I_{|\mathcal{S}|})  \mathbf{v}^t  - \left[ \left(\frac{1}{n}\mathbf{e}\mathbf{e}^T \otimes I_{|\mathcal{S}|}\right) \right] (W^t \otimes I_{|\mathcal{S}|}) \mathbf{v}^t  + \left[I_{n|\mathcal{S}|} - \left(\frac{1}{n}\mathbf{e}\mathbf{e}^T \otimes I_{|\mathcal{S}|}\right) \right] \mathbf{q}^t
\end{align}
where we have used the recursion for the stacked value vector in \eqref{eq:stacked_value_mixing}, and distributed terms. Now, use the fact that the matrix $W^t$ is doubly stochastic and symmetric to note that $\left[I - \frac{1}{n}\mathbf{e}\mathbf{e}^T\right] W^t =W^t \left[I - \frac{1}{n}\mathbf{e}\mathbf{e}^T\right] $, which is unaffected by the Kronecker products with identity, to group terms as 
\begin{align}\label{eq:consensus_error_recursion2}
	\delta^{t+1} 
	&= (W^t \otimes I_{|\mathcal{S}|}) \left[I_{n|\mathcal{S}|} - \left(\frac{1}{n}\mathbf{e}\mathbf{e}^T\otimes I_{|\mathcal{S}|}\right) \right] \mathbf{v}^t  + \left[I_{n|\mathcal{S}|} - \left(\frac{1}{n}\mathbf{e}\mathbf{e}^T \otimes I_{|\mathcal{S}|}\right) \right] \mathbf{q}^t\nonumber \\
	&=  (W^t \otimes I_{|\mathcal{S}|}) \mathbf{\delta}^{t} + \left[I_{n|\mathcal{S}|} - \left(\frac{1}{n}\mathbf{e}\mathbf{e}^T \otimes I_{|\mathcal{S}|}\right) \right] \mathbf{q}^t		%
\end{align}
Next, recursively apply this logic backwards in time to obtain:
\begin{align}\label{eq:consensus_error_recursion3}
	\delta^{t+1} 
	&= \left[\Phi(t,s) \otimes I_{|\mathcal{S}|}\right]\delta^s + \sum_{l=s}^{t-1}\left(\Phi(t,l+1) \otimes I_{|\mathcal{S}|}\right) \left[I_{n|\mathcal{S}|} - \frac{1}{n}\mathbf{e}\mathbf{e}^T\otimes I_{|\mathcal{S}|}\right] \mathbf{q}^l + \left[I_{n|\mathcal{S}|} - \frac{1}{n}\mathbf{e}\mathbf{e}^T \otimes I_{|\mathcal{S}|}\right] \mathbf{q}^t
\end{align}
Now, note that $[(\frac{1}{n}\mathbf{e}\mathbf{e}^T)\otimes I_{|\mathcal{S}|}]\delta^t=\mathbf{0}$
and 
$[(\frac{1}{n}\mathbf{e}\mathbf{e}^T)\otimes I_{|\mathcal{S}|}][ I_{n|\mathcal{S}|} - (\frac{1}{n}\mathbf{e}\mathbf{e}^T)\otimes I_{|\mathcal{S}|}]\mathbf{q}^t=\mathbf{0}$
, where $\mathbf{0}$ is the vector of all $0$'s in $\mathbb{R}^{n|\mathcal{S}|}$ and $\mathbf{e}\in\mathbb{R}^{n}$ is the vector of all $1$'s as defined earlier. These identities may be seen via the definitions of $\delta^t$ and $\mathbf{q}^t$, respectively, in \eqref{eq:consensus_error_vector_identity} and \eqref{eq:auxiliary_sequence}. Substitute these two identities in the first and second terms on the right-hand side of the previous expression to obtain
\begin{align}\label{eq:consensus_error_recursion4}
	\delta^{t+1} 
	&= \left[\left(\Phi(t,s) - \frac{1}{n}\mathbf{e}\mathbf{e}^T\right) \otimes I_{|\mathcal{S}|}\right]\delta^s 
	+ \sum_{l=s}^{t-1}
	\left(\left(\Phi(t,l+1)
	- \frac{1}{n}(\mathbf{e}\mathbf{e}^T) \right)\otimes I_{|\mathcal{S}|} \right)
	\left[I_{n|\mathcal{S}|} - \frac{1}{n}\mathbf{e}\mathbf{e}^T\otimes I_{|\mathcal{S}|}\right] \mathbf{q}^l  \nonumber \\
	&\qquad+ \left[I_{n|\mathcal{S}|} - \frac{1}{n}\mathbf{e}\mathbf{e}^T \otimes I_{|\mathcal{S}|}\right] \mathbf{q}^t
\end{align}
where $\Phi(t,s)$ is defined in \eqref{transition_matrix} as the product of weight matrices. Now, compute the norm of both sides, and apply the triangle inequality together with Lemma \ref{tran_bound} to upper-bound the difference of product matrices to the $n\times n$ uniform weight matrix $\frac{1}{n}\mathbf{e}\mathbf{e}^T$ as
\begin{align}\label{eq:consensus_error_recursion5}
	\|\delta^{t+1} \|
	&\leq \left\|\left(\Phi(t,0) - \frac{1}{n}\mathbf{e}\mathbf{e}^T\right) \otimes I_{|\mathcal{S}|}\right\|_F
	\|\delta^0\| + \sum_{l=0}^{t-1}\left\|\left(\Phi(t,l+1)
	- \frac{1}{n}(\mathbf{e}\mathbf{e}^T) \right)\otimes I_{|\mathcal{S}|}\right\|_F \left\|I_{n|\mathcal{S}|}
	- \frac{1}{n}(\mathbf{e}\mathbf{e}^T)\otimes I_{|\mathcal{S}|}\right\|_F \|\mathbf{q}^l\| \nonumber \\
	&\qquad+ \left\|I_{n|\mathcal{S}|} - \frac{1}{n}\mathbf{e}\mathbf{e}^T \otimes I_{|\mathcal{S}|}\right\|_F \|\mathbf{q}^t\| \nonumber \\
	&\leq \Gamma \rho^{t}\|\delta^0\|  + \sum_{l=0}^{t-1}\Gamma \rho^{t-(l+1)} \|\mathbf{q}^l\| + \|\mathbf{q}^t\|
\end{align}
which, after computing the conditional expectation of both sides and applying \eqref{eq:upper_bound_mean_q} yields
\begin{align}\label{eq:consensus_error_recursion5}
	\mathbb{E}\left[\|\delta^{t+1}\| \mid \mathcal{F}_t\right] 
	&\leq \Gamma \rho^{t}\|\delta^0\| +\mathcal{O}(\sqrt{n}\alpha)\left[1+ \sum_{l=0}^{t-1}\Gamma \rho^{t-(l+1)}\right]
\end{align}
%
Now, under initialization $v_i^t=0$ for all $i$, we have $\delta^0=0$. Now, we evaluate the finite geometric sum via the fact that  $0<\rho<1$. In particular,$ \sum_{l=0}^{t-1} \rho^{t-(l+1)} =  \sum_{u=0}^{t-1} \rho^{u}=(1-\rho^{t-1})/(1-\rho)$. Together with \eqref{eq:consensus_error_recursion5}, we may conclude Lemma \ref{dist_mean}.  $\hfill \square$

\subsection{Proof of Lemma \ref{dist_mean} Statement (ii)}\label{proof_dist_mean2}
Begin by employing the definition of the the auxiliary sequences $\{\tilde p_i^t\}_{t\geq 0}$ and $\{\widetilde q_i^t\}_{t\geq 0}$ for all agents $i \in \mathcal{V}$. First, we analyze $\widetilde q_i^t$ and show its magnitude is upper-bounded. In particular, note that
\begin{align}\label{eq:upper_bound_q2}
	\widetilde q_i^t &= \mu_i^{t+1}- \widetilde{\mu}_i^t,
\end{align}
which would then imply that $\mu_i^{t+1}=\widetilde{\mu}_i^t+ \widetilde q_i^t$. {In a key departure from the standard analysis of the consensus error  in the proof of statement (i) of Lemma \ref{dist_mean}, we refine the discrepancy to instead be in terms of KL divergence for the occupancy measure. Specifically, consider the term $\|\mu_i^{t+1}- \widetilde{\mu}_i^t\|_1$. Via Pinsker's inequality, we can write}
\begin{align}\label{here22222222}
	\|\mu_i^{t+1}- \widetilde{\mu}_i^t\|_1^2=\|\widetilde{\mu}_i^t- \mu_i^{t+1}\|_1^2\leq 2D_{KL}(\widetilde{\mu}_i^t||\mu_i^{t+1}).
\end{align}
Since we have $\widetilde{\mu}_i^t$ as the weighted average of ${\mu}_i^t\in\mathcal{U}$, which implies that $\widetilde{\mu}_i^t\in\mathcal{U}$. Therefore, $D_{KL}(\widetilde{\mu}_i^t||\mu_i^{t+1})\leq D_{KL}(\widetilde{\mu}_i^t||\mu_i^{t+\frac{1}{2}})$ where $\mu_i^{t+\frac{1}{2}}\notin\mathcal{U}$ but it is a valid distribution. Hence, we could further write \eqref{here22222222} as 
\begin{align}\label{here3}
	\|\mu_i^{t+1}- \widetilde{\mu}_i^t\|_1^2=\|\widetilde{\mu}_i^t- \mu_i^{t+1}\|_1^2\leq & 2D_{KL}(\widetilde{\mu}_i^t||\mu_i^{t+\frac{1}{2}})\nonumber
	\\
	=&2 \sum_{s\in\mathcal{S}}\sum_{a\in\mathcal{A}}\widetilde\mu^t(s,a)\log\left(\frac{\widetilde\mu^t(s,a)}{\mu_i^{t+\frac{1}{2}}(s,a)}\right)\nonumber\\
	=& 2\sum_{s\in\mathcal{S}}\sum_{a\in\mathcal{A}} \widetilde\mu^t(s,a) \log \left(\frac{Z_i}{\exp(\Delta_i^{t+1}(s,a))} \right)
	\nonumber\\
	=& 2\log(Z_i)-2\sum_{s\in\mathcal{S}}\sum_{a\in\mathcal{A}} \widetilde\mu^t(s,a)  \left(\Delta_i^{t+1}(s,a) \right).
\end{align}
Note the bound on $\log(Z_i)$ from \eqref{eq:log_simplify2} to obtain
\begin{align} \label{eq:log_simplify2222}
	\|\mu_i^{t+1}- \widetilde{\mu}_i^t\|_1^2   
	&\leq 2\sum_{s\in\mathcal{S}}\sum_{a\in\mathcal{A}}  \widetilde{\mu}_i^t(s,a) \Delta_i^{t+1}(s,a) +  \sum_{s\in\mathcal{S}}\sum_{a\in\mathcal{A}}  \widetilde{\mu}_i^t(s,a) (\Delta_i^{t+1}(s,a))^2-2\sum_{s\in\mathcal{S}}\sum_{a\in\mathcal{A}} \widetilde\mu^t(s,a)  \left(\Delta_i^{t+1}(s,a) \right)
	\nonumber
	\\
	&=\sum_{s\in\mathcal{S}}\sum_{a\in\mathcal{A}}  \widetilde{\mu}_i^t(s,a) (\Delta_i^{t+1}(s,a))^2.
\end{align}
From the statement of Lemma 5, we can write 
\begin{align}\label{eq:upper_bound_mean_q22}
	\mathbb{E}[	\|\mu_i^{t+1}- \widetilde{\mu}_i^t\|_1^2   \mid \mathcal{F}_t] = 	&\sum_{s\in\mathcal{S}}\sum_{a\in\mathcal{A}}  \widetilde{\mu}_i^t(s,a) \mathbb{E}[	(\Delta_i^{t+1}(s,a))^2   \mid \mathcal{F}_t]  
	\nonumber
	\\
	\leq &  4(4t_{mix}+1)^2\beta^2|\mathcal{S}| |\mathcal{A}|.
\end{align}
From Jensen's inequality, we note that 
\begin{align}
	\left(\mathbb{E}[	\|\mu_i^{t+1}- \widetilde{\mu}_i^t\|_1   \mid \mathcal{F}_t]\right)^2\leq& 	\mathbb{E}[	\|\mu_i^{t+1}- \widetilde{\mu}_i^t\|_1^2   \mid \mathcal{F}_t]\nonumber
	\\
	\leq& 4(4t_{mix}+1)^2\beta^2|\mathcal{S}| |\mathcal{A}|.
\end{align}
Taking square root on both sides, we get
\begin{align} \label{main_here5}
	\mathbb{E}[	\|\mu_i^{t+1}- \widetilde{\mu}_i^t\|_1   \mid \mathcal{F}_t]
	\leq& 2(4t_{mix}+1)\beta\sqrt{|\mathcal{S}| |\mathcal{A}|}.
\end{align}
In order to proceed in a similar manner to Lemma \ref{dist_mean} Statement (i) for consensus error proof for $\mu$, we need to bound
\begin{align}\label{eq:upper_bound_mean_q}
	\mathbb{E}[\| \widetilde q_i^t\| \mid \mathcal{F}_t] = 	&\mathbb{E}\left[\| \widetilde{\mu}_i^t-\mu_i^{t+1}\|_1~|~\mathcal{F}_t\right]
	\nonumber
	\\
	\leq & \mathcal{O}({\beta}),
\end{align}

where the second inequality holds from \eqref{main_here5}. Note that we may stack $\widetilde q_i^t$ over all $i\in\mathcal{V}$ to write $\mathbb{E}[\|\widetilde{\mathbf{q}}_t\| \mid \mathcal{F}_t] \leq \mathcal{O}(\sqrt{n} \beta)$, where $\widetilde{\mathbf{q}}_t=[[\widetilde q_1^t]^T ;  \cdots ; [\widetilde q_n^t]^T]\in\mathbb{R}^{n|\mathcal{S}||\mathcal{A}|}$ stacks $\widetilde q_i^t\in\mathbb{R}^{|\mathcal{S}||\mathcal{A7}|}$ [cf. \eqref{eq:auxiliary_sequence}] over all $i$. For the moment, consider that the value vector is scalar for each agent. Then, using this definition of $\widetilde q_i^t$, we may write the the distribution vector $\boldsymbol{\mu}^t=[[\mu_1^t]^T ;  \cdots ; [\mu_n^t]^T]\in\mathbb{R}^{n|\mathcal{S}||\mathcal{A}|}$ stacking $\mu_i^t$ across all $i$ as
\begin{align}\label{eq:stacked_value_mixing23}
	\boldsymbol{\mu}^{t+1} &= \widetilde{\boldsymbol{\mu}}^t + \widetilde{\mathbf{q}}^t\nonumber \\
	&= (W^t \otimes I_{|\mathcal{S}||\mathcal{A}|}) \boldsymbol{\mu}^{t}  +\widetilde{\mathbf{q}}^t,
\end{align}
where  $\otimes $ denotes the Kronecker product and $I_{|\mathcal{S}||\mathcal{A}|}\in\mathbb{R}^{|\mathcal{S}||\mathcal{A}|\times |\mathcal{S}||\mathcal{A}|}$ is the $|\mathcal{S}||\mathcal{A}|\times |\mathcal{S}||\mathcal{A}|$ identity matrix. 
With these observations, shift to defining the respective sequences that track the distance to the mean [cf. \eqref{eq:average_vector}]
\begin{equation}\label{eq:consensus_error_vector2}
	\bar\delta_i^t = \mu_i^t - \overline{\mu}^t \; . 
\end{equation}
We first analyze the evolution of $\bar\delta_i^t$ stacked over all $i$ defined as $\boldsymbol{\bar\delta}^t=[[\delta_1^t]^T ;  \cdots ; [\delta_n^t]^T]\in\mathbb{R}^{n|\mathcal{A}|}$.
Next, following the similar steps from \eqref{eq:stacked_value_mixing} to \eqref{eq:consensus_error_recursion5}, we can write

\begin{align}\label{eq:consensus_error_recursion53}
	\mathbb{E}\left[\|\boldsymbol{\bar\delta}^{t+1}\| \mid \mathcal{F}_t\right] 
	&\leq \Gamma \rho^{t}\|\boldsymbol{\bar\delta}^{0}\| +(2(4t_{mix}+1)\beta\sqrt{|\mathcal{S}| |\mathcal{A}|} \sqrt{n})\left[1+ \sum_{l=0}^{t-1}\Gamma \rho^{t-(l+1)}\right]
\end{align}
%
Now, under initialization $u_i^0=\zeta$ for all $i$, we have $\boldsymbol{\bar\delta}^{0}=0$. Now, we evaluate the finite geometric sum via the fact that  $0<\rho<1$. In particular,$ \sum_{l=0}^{t-1} \rho^{t-(l+1)} =  \sum_{u=0}^{t-1} \rho^{u}=(1-\rho^{t-1})/(1-\rho)$. Together with \eqref{eq:consensus_error_recursion53}, we may conclude Lemma \ref{dist_mean} Statement (ii).  $\hfill \square$

\section{Proof of Lemma \ref{lemma:KL_divergence_contraction}}  \label{proof_lemma:KL_divergence_contraction}

%
Begin by expanding the $i$th term on the left-hand side of \eqref{eq:KL_divergence_contraction} using the definition of KL divergence:
\begin{align}\label{eq:KL_divergence_contraction2}
	D_{KL}(\mu^* \| \mu_i^{t+1})-D_{KL}(\mu^* \| \mu_i^{t})  &	\leq 	D_{KL}(\mu^* \| \mu_i^{t+\frac{1}{2}})-D_{KL}(\mu^* \| \mu_i^{t}),
\end{align}
where the inequality holds due to the fact that  $\mu_i^{t+1}$ is the projected version of $\mu_i^{t+\frac{1}{2}}$ onto the space $\mathcal{U}$. Next, after expanding the definition of KL divergence, we obtain
\begin{align}\label{eq:KL_divergence_contraction222}
	D_{KL}(\mu^* \| \mu_i^{t+1})-D_{KL}(\mu^* \| \mu_i^{t})  &	\leq  \sum_{s\in\mathcal{S}}\sum_{a\in\mathcal{A}}\mu^*(s,a)\log\left(\frac{\mu^*(s,a)}{\mu_i^{t+\frac{1}{2}}(s,a)}\right)-\sum_{s\in\mathcal{S}}\sum_{a\in\mathcal{A}}\mu^*(s,a)\log\left(\frac{\mu^*(s,a)}{\mu_i^{t}(s,a)}\right)
	\nonumber
	\\
	&	= \sum_{s\in\mathcal{S}}\sum_{a\in\mathcal{A}}\mu^*(s,a)\log\left(\frac{\mu_i^{t}(s,a)}{\mu_i^{t+\frac{1}{2}}(s,a)}\right).
\end{align}
Add the term $\sum_{s\in\mathcal{S}}\sum_{a\in\mathcal{A}}\mu^*(s,a)\log\left(\frac{\widetilde\mu^t(s,a)}{\widetilde\mu^t(s,a)}\right)$ in \eqref{eq:KL_divergence_contraction222} to obtain
\begin{align}\label{eq:KL_divergence_contraction3}
	D_{KL}(\mu^* \| \mu_i^{t+1})-D_{KL}(\mu^* \| \mu_i^{t})  &	\leq \sum_{s\in\mathcal{S}}\sum_{a\in\mathcal{A}}\mu^*(s,a)\log\left(\frac{\widetilde\mu^t(s,a)}{\mu_i^{t+\frac{1}{2}}(s,a)}\right)+\sum_{s\in\mathcal{S}}\sum_{a\in\mathcal{A}}\mu^*(s,a)\log\left(\frac{\mu_i^{t}(s,a)}{\widetilde\mu^t(s,a)}\right).
\end{align}
Take the average over $i$ to obtain
\begin{align}\label{eq:KL_divergence_contraction43}
	\frac{1}{n}\sum_{i=1}^nD_{KL}(\mu^* \| \mu_i^{t+1})-	\frac{1}{n}\sum_{i=1}^nD_{KL}(\mu^* \| \mu_i^{t})  	\leq & 	\frac{1}{n}\sum_{i=1}^n \sum_{s\in\mathcal{S}}\sum_{a\in\mathcal{A}}\mu^*(s,a)\log\left(\frac{\widetilde\mu^t(s,a)}{\mu_i^{t+\frac{1}{2}}(s,a)}\right)
	\nonumber
	\\
	&+	\underbrace{\frac{1}{n}\sum_{i=1}^n\sum_{s\in\mathcal{S}}\sum_{a\in\mathcal{A}}\mu^*(s,a)\log\left(\frac{\mu_i^{t}(s,a)}{\widetilde\mu^t(s,a)}\right)}_{T_0}.
\end{align}
Let us consider the second term on the right hand side of \eqref{eq:KL_divergence_contraction43} as follows:
\begin{align}
	T_0= \frac{1}{n}\sum_{i=1}^n\sum_{s\in\mathcal{S}}\sum_{a\in\mathcal{A}}\mu^*(s,a)\log\left({\mu_i^{t}(s,a)}\right)-\frac{1}{n}\sum_{i=1}^n\sum_{s\in\mathcal{S}}\sum_{a\in\mathcal{A}}\mu^*(s,a)\log\left({\widetilde\mu^t(s,a)}\right).
\end{align}
Substitute the expression for $\widetilde\mu^t(s,a)$ (cf. \eqref{eq:consensus_round}) to obtain
\begin{align}
	T_0= \frac{1}{n}\sum_{i=1}^n\sum_{s\in\mathcal{S}}\sum_{a\in\mathcal{A}}\mu^*(s,a)\log\left({\mu_i^{t}(s,a)}\right)-\frac{1}{n}\sum_{i=1}^n\sum_{s\in\mathcal{S}}\sum_{a\in\mathcal{A}}\mu^*(s,a)\log\left( \sum_{j=1}^{n} w_{ij}^t \mu_j^t(s,a)\right).
\end{align}
Since $-\log(\cdot)$ is convex, via Jensen's inequality, we obtain
\begin{align}
	T_0\leq  \frac{1}{n}\sum_{i=1}^n\sum_{s\in\mathcal{S}}\sum_{a\in\mathcal{A}}\mu^*(s,a)\log\left({\mu_i^{t}(s,a)}\right)-\frac{1}{n}\sum_{i=1}^n\sum_{s\in\mathcal{S}}\sum_{a\in\mathcal{A}}\mu^*(s,a)\sum_{j=1}^{n} w_{ij}^t \log\left( \mu_j^t(s,a)\right).
\end{align}
After changing the order of summation, we obtain
\begin{align}
	T_0\leq  \frac{1}{n}\sum_{i=1}^n\sum_{s\in\mathcal{S}}\sum_{a\in\mathcal{A}}\mu^*(s,a)\log\left({\mu_i^{t}(s,a)}\right)-\frac{1}{n}\sum_{j=1}^n\underbrace{\sum_{i=1}^{n} w_{ij}^t}_{=1}\sum_{s\in\mathcal{S}}\sum_{a\in\mathcal{A}}\mu^*(s,a) \log\left( \mu_j^t(s,a)\right).
\end{align}
After simplifications, we obtain
\begin{align}
	T_0\leq & \frac{1}{n}\sum_{i=1}^n\sum_{s\in\mathcal{S}}\sum_{a\in\mathcal{A}}\mu^*(s,a)\log\left({\mu_i^{t}(s,a)}\right)-\frac{1}{n}\sum_{j=1}^n\sum_{s\in\mathcal{S}}\sum_{a\in\mathcal{A}}\mu^*(s,a) \log\left( \mu_j^t(s,a)\right)
	\nonumber
	\\
	=& 0.\label{here2222}
\end{align}
Using the upper bound in \eqref{here2222} into \eqref{eq:KL_divergence_contraction43} to obtain
\begin{align}\label{eq:KL_divergence_contraction433}
	\frac{1}{n}\sum_{i=1}^nD_{KL}(\mu^* \| \mu_i^{t+1})-\frac{1}{n}\sum_{i=1}^nD_{KL}(\mu^* \| \mu_i^{t})  &	\leq \underbrace{\frac{1}{n}\sum_{i=1}^n \sum_{s\in\mathcal{S}}\sum_{a\in\mathcal{A}}\mu^*(s,a)\log\left(\frac{\widetilde\mu^t(s,a)}{\mu_i^{t+\frac{1}{2}}(s,a)}\right)}_{:=T_1}.
\end{align}
Let us consider the term $T_1$, which is essentially similar to the quantity that appears in \citep{Wang2020}[Lemma A.2]. To upper-bound this estimate, make use of the definition of $\mu_i^{t + \frac{1}{2}}(s,a) = \frac{\widetilde{\mu}_i^t(s,a) \exp( \Delta_i^{t+1}(s,a))}{\sum_{s'} \sum_{a'} \widetilde{\mu}_i^t(s,a) \exp( \Delta_i^{t+1}(s',a'))}$ in Algorithm \ref{alg_single} as follows:
\begin{align}\label{eq:define_Z}
	T_1  = &\frac{1}{n}\sum_{i=1}^n \sum_{s\in\mathcal{S}}\sum_{a\in\mathcal{A}}\mu^*(s,a)\log\left(\frac{\widetilde\mu^t(s,a)}{\mu_i^{t+\frac{1}{2}}(s,a)}\right)\nonumber\\
	=& \frac{1}{n}\sum_{i=1}^{n}\sum_{s\in\mathcal{S}}\sum_{a\in\mathcal{A}} \mu^*(s,a) \log \left(\frac{Z_i}{\exp(\Delta_i^{t+1}(s,a))} \right)
\end{align}
where we have defined $Z_i $$=$$ \sum_{s \in S} \sum_{a \in A} \widetilde{\mu}_i^t(s,a) \exp(\Delta_i^{t+1}(s,a))$, and used the fact that a log of a product equals the sum of logs. Next, use the fact that a log of a ratio is the difference of the logs in the last term on the right-hand side of \eqref{eq:define_Z} to write:
\begin{align}\label{eq:lemma_primal_dual_proof_start22}
	T_1    &= \frac{1}{n}\sum_{i=1}^{n}\sum_{s\in\mathcal{S}}\sum_{a\in\mathcal{A}}    \mu^*(s,a) \log (Z_i) -\frac{1}{n}
	\sum_{i=1}^{n}\sum_{s\in\mathcal{S}}\sum_{a\in\mathcal{A}}   \mu^*(s,a) \Delta_i^{t+1}(s,a) \nonumber\\
	&= \frac{1}{n}\sum_{i=1}^{n} \log (Z_i) - \frac{1}{n}\sum_{i=1}^{n}\sum_{s\in\mathcal{S}}\sum_{a\in\mathcal{A}}  \mu^*(s,a) \Delta_i^{t+1}(s,a)
\end{align}
where we have used the fact that $\mu^*(s,a)$ is a likelihood whose sum over all $s\in\mathcal{S},a\in\mathcal{A}$ is unit to simplify the third term on the right-hand side. 
Next, to analyze  $Z_i$, especially its logarithm, it turns out to be useful to establish that  $\Delta_i^{t+1}(s,a)$ [cf. \eqref{eq:dual_gradient}] is always either negative or zero. To do so, note that $v_i^t(s) \in [-2t^*_{mix}, 2t^*_{mix}]$, and $r_i^t(s,s',a) \in [0, 1]$, and $M $$=$$ 4t_{mix}+1$ Thus, $v_i^t(s')-v_i^t(s)+r_i^t(s,s',a)-M $$\leq$$ 2t^*_{mix} - (-2t^*_{mix}) +1 -(4t^*_{mix} +1) $$\leq$$ 0$. Then it follows that $\Delta_i^{t+1}(s,a) \leq 0$ for all $s \in S$ and $a \in A$ with probability 1. 
Now, we shift to considering the definition of $Z_i$, whose logarithm in \eqref{eq:lemma_primal_dual_proof_start22} is given as:
\begin{align}\label{eq:log_simplify1}
	\log (Z_i) &= \log\left(\sum_{s\in\mathcal{S}}\sum_{a\in\mathcal{A}}  \widetilde{\mu}_i^t(s,a) \exp(\Delta_i^{t+1}(s,a))\right) \\
	&\leq \log \left(\sum_{s\in\mathcal{S}}\sum_{a\in\mathcal{A}}  \widetilde{\mu}_i^t(s,a) (1+ \Delta_i^{t+1}(s,a) + \frac{1}{2} (\Delta_i^{t+1}(s,a))^2 )\right) \nonumber
\end{align}
where we have used the inequality $e^{x}\leq 1+x+\frac{1}{2}x^2$ for $x\leq 0$ above. Next, let us apply $\log(1+x) \leq x$ which holds for all $x\in\mathbb{R}$ to the right-hand side as:
\begin{align} \label{eq:log_simplify2}
	\log Z_i    
	&\leq \log\left(1+ \sum_{s\in\mathcal{S}}\sum_{a\in\mathcal{A}}  \widetilde{\mu}_i^t(s,a) \Delta_i^{t+1}(s,a) + \frac{1}{2} \sum_{s\in\mathcal{S}}\sum_{a\in\mathcal{A}}  \widetilde{\mu}_i^t(s,a) (\Delta_i^{t+1}(s,a))^2\right) \nonumber \\
	&\leq \sum_{s\in\mathcal{S}}\sum_{a\in\mathcal{A}}  \widetilde{\mu}_i^t(s,a) \Delta_i^{t+1}(s,a) + \frac{1}{2} \sum_{s\in\mathcal{S}}\sum_{a\in\mathcal{A}}  \widetilde{\mu}_i^t(s,a) (\Delta_i^{t+1}(s,a))^2.
\end{align}
Via the analysis of $\log (Z_i)$ in \eqref{eq:log_simplify2},  $T_1$ defined in \eqref{eq:lemma_primal_dual_proof_start22}, , simplifies to
\begin{align}\label{eq:lemma_primal_dual_proof_start222}
	T_1   
	&\leq  \frac{1}{n}\sum_{i=1}^{n} \sum_{s\in\mathcal{S}}\sum_{a\in\mathcal{A}}  \widetilde{\mu}_i^t(s,a) \Delta_i^{t+1}(s,a) + \frac{1}{2n}\sum_{i=1}^{n} \sum_{s\in\mathcal{S}}\sum_{a\in\mathcal{A}}  \widetilde{\mu}_i^t(s,a) (\Delta_i^{t+1}(s,a))^2 - \frac{1}{n}\sum_{i=1}^{n}\sum_{s\in\mathcal{S}}\sum_{a\in\mathcal{A}}  \mu^*(s,a) \Delta_i^{t+1}(s,a)
	\nonumber
	\\
	&= \frac{1}{n} \sum_{i=1}^{n} \sum_{s\in\mathcal{S}}\sum_{a\in\mathcal{A}}\left(\widetilde{\mu}_i^t(s,a)-{\mu}^*(s,a)\right) \Delta_i^{t+1}(s,a) + \frac{1}{2n}\sum_{i=1}^{n} \sum_{s\in\mathcal{S}}\sum_{a\in\mathcal{A}}  \widetilde{\mu}_i^t(s,a) (\Delta_i^{t+1}(s,a))^2 .
\end{align}
Let us utilize this bound into the right hand side of \eqref{eq:KL_divergence_contraction433} to conclude
\begin{align}
	\frac{1}{n}\sum_{i=1}^nD_{KL}(\mu^* \| \mu_i^{t+1})-\frac{1}{n}\sum_{i=1}^nD_{KL}(\mu^* \| \mu_i^{t})  	\leq & \frac{1}{n}\sum_{i=1}^{n} \sum_{s\in\mathcal{S}}\sum_{a\in\mathcal{A}}\left(\widetilde{\mu}_i^t(s,a)-{\mu}^*(s,a)\right) \Delta_i^{t+1}(s,a) 
	\nonumber
	\\
	&+ \frac{1}{2}\sum_{i=1}^{n} \sum_{s\in\mathcal{S}}\sum_{a\in\mathcal{A}}  \widetilde{\mu}_i^t(s,a) (\Delta_i^{t+1}(s,a))^2.
\end{align}

as stated in Lemma \ref{lemma:KL_divergence_contraction}. $\hfill \square$

\section{Proof of Lemma \ref{lemma:occupancy_decrement}}\label{proof_lemma:occupancy_decrement}

Consider the update direction of the local occupancy $\Delta_i^{t+1}$ [cf. \eqref{eq:dual_gradient}] in conditional expectation (where we use the fact that the update occurs with probability $\widetilde{\mu}_i^t(s,a)$), stated as:
\begin{align}
	\mathbf{E}[\Delta_i^{t+1}(s,a)| \mathcal{F}_t] &= \beta\widetilde{\mu}_i^t(s,a)\sum_{s' \in S} p_{ss'}(a) \frac{v_i^t(s') - v_i^t(s)}{\widetilde{\mu}_i^t(s,a)} + \beta\widetilde{\mu}_i^t(s,a)\sum_{s' \in S} p_{ss'}(a) \frac{r_i^t(s,s',a) - M}{\widetilde{\mu}_i^t(s,a)}\\
	&=  \beta\left((P_a v_i^t - v_i^t +r_a)(s) - M \right)\label{conditional_Expectation}
\end{align}
where we have canceled out a factor of $\widetilde{\mu}_i^t(s,a)$.
Then compute the quantity in the left-hand side of \eqref{eq:occupancy_decrement} and exploit the fact that that the sum of the elements of occupancy measures is unit, i.e., $\sum_{s\in\mathcal{S},a\in\mathcal{A}}\widetilde{\mu}_i^t(s,a)=\sum_{s\in\mathcal{S},a\in\mathcal{A}}{\mu}_i^t(s,a)=1 $ and $\sum_{a \in A} (\widetilde{\mu}_i^t(a)=\sum_{a \in A} ({\mu}_i^t(a)=1$.
\begin{align*}
	\frac{1}{n}\sum_{i=1}^{n} \sum_{s \in S} \sum_{a \in A} (\widetilde{\mu}_i^t(s,a) - \mu^*(s,a)) \mathbb{E}[\Delta_i^{t+1}(s,a)\mid \mathcal{F}_t] 
	&=  \frac{\beta}{n}\sum_{i=1}^{n} \sum_{a \in A} (\widetilde{\mu}_i^t(a) - \mu^*(a))^T[(P_a - I) v_i^t +r_a)]\nonumber \\
	&= \frac{\beta}{n} \sum_{i=1}^{n}\sum_{s \in S} \sum_{a \in A} (\widetilde{\mu}_i^t(s,a) - \mu^*(s,a))[(P_a v_i^t - v_i^t +r_a)(s) - M] .
\end{align*}
The result is \eqref{eq:occupancy_decrement}.
\hfill $\square$

\section{Proof of Lemma \ref{lemma:mean_square_dual}}\label{proof_lemma_mean_Swuare_dual}
 Consider the expression on the left-hand side of \eqref{eq:mean_square}. Use the definition of $\Delta_i^{t+1}(s,a)$ [cf. \eqref{eq:dual_gradient}] as well as the fact that it occurs with probability  $\widetilde{\mu}_i^t(s,a)$ to write:
\begin{align}\label{eq:lemma_mean_square_start}
	\sum_{s \in S} \sum_{a \in A} \widetilde{\mu}_i^t(s,a) \mathbb{E}[(\Delta_i^{t+1}(s,a))^2\mid\mathcal{F}_t] &= \sum_{s \in S}\sum_{a \in A} \widetilde{\mu}_i^t(s,a) \widetilde{\mu}_i^t(s,a) \sum_{s' \in S} p_{ss'}(a) \Big( \beta \frac{h_i^t(s') - h_i^t(s)+r_i^t(s,s',a)-M}{\widetilde{\mu}_i^t(s,a)}\Big)^2
\end{align}    
Then, we may group the sum together, and cancel out a factor of $\widetilde{\mu}_i^t(s,a)^2$ to write:
\begin{align}    \label{eq:lemma_mean_square_near_final}
	\sum_{s \in S} \sum_{a \in A} \widetilde{\mu}_i^t(s,a) \mathbb{E}[(\Delta_i^{t+1}(s,a))^2\mid\mathcal{F}_t] 
	&= \sum_{s \in S}\sum_{a \in A} \sum_{s' \in S} p_{ss'}(a) (\beta (v_i^t(s') - v_i^t(s)+r_i^t(s,s',a)-M))^2 \nonumber\\
	&\leq \sum_{s \in S}\sum_{a \in A} \sum_{s' \in S} p_{ss'}(a) (2 \beta (4t_{mix}+1))^2\nonumber \\
	&= 4(4t_{mix}+1)^2\beta^2  |\mathcal{S}| |\mathcal{A}| 
\end{align}
where the inequality uses the fact that $v_i^t(s) \in [-2t^*_{mix}, 2t^*_{mix}]$, and $r_i^t(s,s',a) \in [0, 1]$, and $M $$=$$ 4t_{mix}+1$, and the last equality evaluates the sum, using the definition of the transition probability and the cardinality of the spaces, and yields Lemma \ref{lemma:mean_square_dual}.

  \hfill $\square$

Next we establish a decrement-like property on the local differential value vector $v_i^{t}$ with respect to the locally optimal $v^{*}$. 

\section{Proof of Lemma \ref{lemma:value_decrement}} \label{proof_lemma:value_decrement}

 This result is a generalization of the proof of \citep{Wang2020}[Lemma 7] that contains additional consensus-error terms. First, recall the definition of $v_i^{t+1}$ in \eqref{eq:value_updates}. Begin by considering the norm-difference of the global average of this quantity $\overline{v}^t$ [cf. \eqref{eq:average_vector}] to the global optimizer $v^*$ defined by the one that minimizes $\sum_i L_i(\mu, v)$ in \eqref{eq:local_lagrangian}:
%
%
\begin{align}\label{eq:value_decrement_consensus}
	\left\| \frac{1}{n}\sum_{i=1}^n v_i^{t+1} - v^*\right\|^2 &=
	\left\| \frac{1}{n}\sum_{i=1}^n \Pi_{\mathcal{V}} [\widetilde{v}_i^t +  d_i^{t+1}] - v^*\right\|^2 
	\nonumber \\
	&=
	\left\| \frac{1}{n}\sum_{i=1}^n \left(\widetilde{v}_i^t+\Pi_{\mathcal{V}} [\widetilde{v}_i^t +  d_i^{t+1}]-\widetilde{v}_i^t\right) - v^*\right\|^2 \nonumber
	\\
	&= \left\| \frac{1}{n}\sum_{i=1}^n \left(\widetilde{v}_i^t +  q_i^{t+1}\right) - v^*\right\|^2
	%
	%
\end{align}
where $q_i^{t+1}:=\Pi_{\mathcal{V}} [\widetilde{v}_i^t +  d_i^{t+1}]-\widetilde{v}_i^t$ and we make use of the non-expansiveness of the projection.
Continue by making use of the definition of the consensus round $\widetilde{v}_i^t=\sum_{j=1}^{n} w_{ij}^t v_j^t$ in \eqref{eq:consensus_round} in the right-hand side of the preceding expression to write
\begin{align}\label{eq:value_decrement_consensus}
	\left\| \frac{1}{n}\sum_{i=1}^n \sum_{j=1}^n w_{ij}^t {v}_j^t + \frac{1}{n} \sum_{i=1}^n q_i^{t+1} - v^*\right\|^2
	&= \left\| \frac{1}{n}\sum_{j=1}^n \left( \sum_{i=1}^n w_{ij}^t \right) {v}_j^t + \frac{1}{n} \sum_{i=1}^n q_i^{t+1} - v^*\right\|^2 \nonumber \\
	& = \left\|  \overline{v}^t +  \frac{1}{n}\sum_{i=1}^n q_i^{t+1} - v^*\right\|^2 
\end{align}
where we have applied the fact that the mixing matrix $W^t$ is doubly stochastic for all $t$.
Now, let us expand the square on the right-hand side of the previous expression, and make use of the short-hand $\overline{v}^t=(1/n)\sum_j v_j^t$ to write:
\begin{align}\label{eq:value_decrement_mixing}
	\left\|  \overline{v}^t +  \frac{1}{n}\sum_{i=1}^n q_i^{t+1} - v^*\right\|^2 & \leq
	\left\| \overline{v}^t - v^* \right \|^2 + \frac{2 \alpha}{n} \left(\overline{v}^t - v^* \right)^T\left( \frac{1}{n}\sum_{i=1}^n \sum_{a\in\mathcal{A}} (I - P_a)^T\widetilde{\mu}_i^t(a)\right) + \frac{1}{n^2}\sum_{i=1}^n\left\| q_i^{t+1}\right\|^2
\end{align}
Let us consider the term now $\left\| q_i^{t+1}\right\|$ and try to develop the upper bound
\begin{align}
	\left\| q_i^{t+1}\right\|^2=&	\left\|\Pi_{\mathcal{V}} [\widetilde{v}_i^t +  d_i^{t+1}]-\widetilde{v}_i^t\right\|^2
	\\
	\leq & \left\| [\widetilde{v}_i^t +  d_i^{t+1}]-\widetilde{v}_i^t\right\|^2 = \left\| d_i^{t+1}\right\|^2
\end{align}
Now, compute the expectation conditional on filtration $\mathcal{F}_t$ of the previous expression and apply \eqref{eq:mean_square}, we get
\begin{align}
	\mathbb{E}\left[\left\| q_i^{t+1}\mid \mathcal{F}_t\right\|^2\right]\leq & \mathbb{E}\left[\left\| d_i^{t+1}\right\|^2\mid \mathcal{F}_t\right]= \mathcal{O}(\alpha^2).
\end{align}

Next, using the definition of the average vector \eqref{eq:average_vector} to write:
\begin{align}\label{eq:value_decrement_mixing2}
	\mathbb{E}\left[\left\| \overline{v}^{t+1}- v^*\right\|^2 \mid \mathcal{F}_t\right]& \leq
	\left\| \overline{v}^t - v^* \right \|^2 + \frac{2 \alpha}{n} \left(\overline{v}^t - v^* \right)^T\left( \frac{1}{n}\sum_{i=1}^n \sum_{a\in\mathcal{A}} (I - P_a)^T\widetilde{\mu}_i^t(a)\right)  + \frac{1}{n} \mathcal{O}(\alpha^2)
\end{align}
%
%
which is as stated in Lemma \ref{lemma:value_decrement}. \hfill $\square$

\section{Proof of Theorem  \ref{eq:KL_divergence_martingale_final400}} \label{proof_eq:KL_divergence_martingale_final400}
Let us start with the statement of Lemma \ref{lemma:KL_divergence_martingale}, 
\begin{align}\label{eq:KL_divergence_martingale_final000}
	\mathbb{E}\left[\mathcal{E}_{t+1}\mid \mathcal{F}_t \right]
	\leq& \mathcal{E}_{t} -\beta\mathcal{D}_t +\beta^2\widetilde{\mathcal{O}}\left(|\mathcal{S}||\mathcal{A}|t_{mix}^2\right) + \underbrace{\frac{\beta}{n}\sum_{i=1}^{n}\sum_{a \in A}\left[ \left(\overline{v}^t -v_i^t \right)^T\left( (I - P_a)^T(\widetilde{\mu}_i^t(a)- \mu^*(a))\right)\right]}_{I}  
	\nonumber
	\\
	&+ \underbrace{\frac{\beta}{n}\sum_{i=1}^{n}\sum_{a \in A}\left[[(P_a - I) v^* +r_a]^T\left(\widetilde{\mu}_i^t(a)-\overline{\mu}^t(a)\right) \right]}_{II}.
\end{align}
Let us consider the terms $I$ and $II$ and try to derive the upper bounds on them separately. Let us write $I$ as
\begin{align}
	I\leq |I|\leq& \frac{\beta}{n}\sum_{i=1}^{n}\sum_{a \in A}\left[ \|\overline{v}^t -v_i^t\|_{\infty}\cdot\| (I - P_a)^T(\widetilde{\mu}_i^t(a)- \mu^*(a))\|_1\right]
	\nonumber
	\\
	\leq& \frac{\beta}{n}\sum_{i=1}^{n}\sum_{a \in A}\left[ \|\overline{v}^t -v_i^t\|_{\infty}\cdot\| I - P_a^T\|_1\cdot\|(\widetilde{\mu}_i^t(a)- \mu^*(a))\|_1\right],
\end{align}
{which follows from the inequality that $\| (I - P_a)^T(\widetilde{\mu}_i^t(a)- \mu^*(a))\|_1\leq \| I - P_a^T\|_1\cdot\|(\widetilde{\mu}_i^t(a)- \mu^*(a))\|_1$. Note that $\| I - P_a^T\|_1\leq \| I \|_1+\| P_a^T\|_1$ and $\| I \|_1=1$ and $\| P_a^T\|_1= \max_{i}\sum_{j=1}^{|S|}| P_a(i,j)|= 1$, hence we could write  }
\begin{align}
	I\leq& \frac{2\beta}{n}\sum_{i=1}^{n}\sum_{a \in A}\left[ \|\overline{v}^t -v_i^t\|_{\infty}\cdot\|\widetilde{\mu}_i^t(a)- \mu^*(a)\|_1\right] .
\end{align}
{Next note that we have $\|\widetilde{\mu}_i^t(a)- \mu^*(a)\|_1\leq \|\widetilde{\mu}_i^t(a)\|_1+\| \mu^*(a)\|_1$ and since $\mu_i^t(a)$ and $\mu^*(a)$ are marginal probability distributions, we have $\|\widetilde{\mu}_i^t(a)- \mu^*(a)\|_1\leq 2$. Hence, we could write}
\begin{align}
	I\leq& \frac{2\beta}{n}\sum_{i=1}^{n}\|\overline{v}^t -v_i^t\|_{1} \sum_{a \in A}\left[ \|\widetilde{\mu}_i^t(a)\|_\infty+\| \mu^*(a)\|_\infty\right] 
	\nonumber
	\\
	\leq& (4\beta|\mathcal{A}|)\cdot\frac{1}{n}\sum_{i=1}^{n}\|\overline{v}^t -v_i^t\|_{1} 
	\nonumber
	\\
	=& (4\beta|\mathcal{A}|)\cdot\frac{1}{n}\sum_{i=1}^{n}\|\delta_i^t\|_{1} 	\nonumber
	\\
	\leq & (4\beta|\mathcal{A}|)\cdot \|\delta^t\|_{1}.
\end{align}
From the statement of Lemma \ref{dist_mean} and the value of $\alpha=|\mathcal{S}|t^2_{mix}\beta$, it holds that 
\begin{align}
	\mathbb{E}\left[I \mid \mathcal{F}_t\right] \leq& (4\beta|\mathcal{A}|)\mathbb{E}\left[ \|\delta^t\|_{1} \mid \mathcal{F}_t\right]
	\nonumber
	\\
	\leq &	8\beta^2(t^2_{mix}\cdot|\mathcal{S}||\mathcal{A}|)\cdot (\sqrt{n})\left[1+  \frac{\Gamma(1-\rho^{t-1})}{1-\rho}\right].
\end{align}
Let us consider the term $II$ term in \eqref{eq:KL_divergence_martingale_final000} and using Cauchy-Schwartz inequality, we can write 
\begin{align}
	II\leq |II| \leq &	\frac{\beta}{n}\sum_{i=1}^{n}\sum_{a \in A}\|(P_a - I) v^* +r_a\|_{\infty}\cdot \|\overline{\mu}^t(a)-\widetilde{\mu}_i^t(a)\|_1\label{here_2}
	\end{align}
{Next, from the triangle inequality, we can write $\|(P_a - I) v^* +r_a\|_{\infty}\leq \|(P_a - I) v^*\|_{\infty}+\|r_a\|_{\infty}$. From the definition of rewards, we note that $\|r_a\|_{\infty}=1$. Also it holds that $\|(P_a - I) v^*\|_{\infty}\leq \|P_a - I\|_{\infty}\cdot\| v^*\|_{\infty}\leq 2 t_{mix}\|P_a - I\|_{\infty} $ from the definition of $\mathcal{V}$. Further from the triangle inequality, we have $\|P_a - I\|_{\infty}\leq \|P_a\|_{\infty} + \|I\|_{\infty}= \|P_a\|_{\infty} + 1$. Next, we use the definition of matrix $\infty$ norm (maximum row sum) as $\|P_a\|_{\infty}= \max_{i}\sum_{j=1}^{|S|}| P_a(i,j)|= 1$. This would imply that $\|P_a - I\|_{\infty}\leq 2$. Combining all these inequalities, we could write
\begin{align}
	\|(P_a - I) v^* +r_a\|_{\infty}\leq & \|(P_a - I) v^*\|_{\infty}+\|r_a\|_{\infty}\nonumber
	\\
	\leq &\|(P_a - I) v^*\|_{\infty}+1
	\nonumber
	\\
	\leq & \|P_a - I\|_{\infty}\cdot\| v^*\|_{\infty}+1	\nonumber
	\\
	\leq &2 t_{mix}\|P_a - I\|_{\infty} +1
	\nonumber
	\\
	\leq &4 t_{mix} +1. \label{new_latest}
\end{align} 
Using the upper bound in \eqref{new_latest} into \eqref{here_2}, we get 
\begin{align}
	II\leq |II| \leq &	\frac{\beta (4 t_{mix} +1)}{n}\sum_{i=1}^{n}\sum_{a \in A} \|\overline{\mu}^t(a)-\widetilde{\mu}_i^t(a)\|_1
	\nonumber
	\\
	\leq &	\frac{\beta (4 t_{mix} +1)}{n}\sum_{i=1}^{n}\sum_{a \in A} \|\overline{\mu}^t(a)-\sum_{j=1}^{n}w_{ij}^t{\mu}_i^t(a)\|_1
	\nonumber
	\\
	\leq &	\beta (4 t_{mix} +1)\sum_{a \in A}\frac{1}{n}\sum_{i=1}^{n}\|\overline{\mu}^t(a)-{\mu}_i^t(a)\|_1
	\nonumber
	\\
	\leq &\beta (4 t_{mix} +1)\frac{1}{n}\sum_{i=1}^{n}\|\overline{\mu}^t-{\mu}_i^t\|_1
	\nonumber
	\\
	\leq &\beta \tilde{t}_{mix}\frac{1}{n}\sum_{i=1}^{n}\|\overline{\mu}^t-{\mu}_i^t\|_1,
\end{align}
where $\tilde{t}_{mix}:=4 t_{mix} +1$. }
%
where $\overline{\mu}^t$ and ${\mu}_i^t$ concatenated the values for all $a\in\mathcal{A}$. To proceed next, we invoke the statement (ii) of Lemma \ref{dist_mean} for $\mu$ to obtain
\begin{align}
	\mathbb{E}\left[II \mid \mathcal{F}_t\right] \leq &	{\beta {\tilde{t}_{mix}}}\cdot (\sqrt{n}2(\tilde{t}_{mix})\beta\sqrt{|\mathcal{S}| |\mathcal{A}|})\cdot\left[1+  \frac{\Gamma(1-\rho^{t-1})}{1-\rho}\right]
	\nonumber
	\\
	\leq&	2\beta^2\left(\tilde{t}_{mix}^2\cdot\sqrt{|\mathcal{A}||\mathcal{S}|}\right)\cdot (\sqrt{n})\cdot\left[1+  \frac{\Gamma(1-\rho^{t-1})}{1-\rho}\right].\label{here_4}
\end{align}
Hence, finally from \eqref{here_2} and \eqref{here_4} combining with the fact that $t^2_{mix}\leq \tilde{t}^2_{mix}$,  we can write the expression in \eqref{eq:KL_divergence_martingale_final000} as 
\begin{align}\label{eq:KL_divergence_martingale_final2}
	\mathbb{E}\left[\mathcal{E}_{t+1}\mid \mathcal{F}_t \right]
	\leq& \mathcal{E}_{t} -\beta\mathcal{D}_t +\beta^2\widetilde{\mathcal{O}}\left(|\mathcal{S}||\mathcal{A}|t_{mix}^2\right) 	\nonumber
	\\
	&+\beta^2\tilde{t}_{mix}^2\left(8|\mathcal{S}||\mathcal{A}|+2\sqrt{|\mathcal{A}||\mathcal{S}|}\right)\cdot (\sqrt{n})\cdot\left[1+  \frac{\Gamma(1-\rho^{t-1})}{1-\rho}\right]
	\nonumber
	\\
	\leq& \mathcal{E}_{t} -\beta\mathcal{D}_t +\beta^2\left(\widetilde{\mathcal{O}}\left(|\mathcal{S}||\mathcal{A}|{t}_{mix}^2\right) +\widetilde{\mathcal{O}}\left(\sqrt{n} |\mathcal{S}||\mathcal{A}| \tilde{t}^2_{mix}\right)\cdot \big[1+  \frac{\Gamma(1-\rho^{t-1})}{1-\rho}\big]\right).
\end{align}
After rearrangement and taking summation over $t$, we get 
\begin{align}\label{eq:KL_divergence_martingale_final3}
	\sum_{t=0}^{T-1}\beta\mathcal{D}_t		\leq& \mathcal{E}_{0}  +\beta^2T\left(\widetilde{\mathcal{O}}\left(|\mathcal{S}||\mathcal{A}|{t}_{mix}^2\right) +\widetilde{\mathcal{O}}\left(\sqrt{n} |\mathcal{S}||\mathcal{A}| \tilde{t}^2_{mix}\right)\cdot \big[  \frac{1+\Gamma}{1-\rho}\big]\right),
\end{align}
where we used that $(1-\rho^{t-1})\leq 1$. Let us define $D(\Gamma, \rho):=\big[  \frac{1+\Gamma}{1-\rho}\big]$ and divide both sides by $\beta T$ to get
\begin{align}\label{eq:KL_divergence_martingale_final4}
	\frac{1}{T}\sum_{t=0}^{T-1}\mathcal{D}_t		\leq& \frac{\mathcal{E}_{0} }{\beta T} +\beta\left(\widetilde{\mathcal{O}}\left({\sqrt{n}} |\mathcal{S}||\mathcal{A}| \tilde{t}^2_{mix}D(\Gamma, \rho)\right)\right).
\end{align}
By selecting the optimal $\beta=\mathcal{\widetilde{\mathcal{O}}}\left(\sqrt{\frac{\mathcal{E}_{0}}{{ {n} |\mathcal{S}||\mathcal{A}| \tilde{t}^2_{mix}D(\Gamma, \rho)}T}}\right)$, we obtain
\begin{align}\label{eq:KL_divergence_martingale_final42}
\mathbb{E}\left[	\frac{1}{T}\sum_{t=0}^{T-1}\mathcal{D}_t\right]		\leq& \widetilde{\mathcal{O}}\left( \tilde{t}_{mix}\sqrt{\frac{\sqrt{n}\mathcal{E}_{0} |\mathcal{S}||\mathcal{A}|D(\Gamma, \rho)}{T}}\right).
\end{align}

\section{Proof of Lemma \ref{prop}} \label{prob_1}

	Let us denote the running average of policy as $\hat\pi$=$\frac{1}{T}\sum_{t=1}^{T}\overline\pi^t$ and note that $     \frac{1}{\sqrt{\tau}|\mathcal{S}|}\mathbf{e} \leq\xi^{\hat\pi}\leq \frac{\sqrt{\tau}}{|\mathcal{S}|}\mathbf{e} $. Hence we can write
	\begin{align}
		\xi^{\hat\pi}\leq \frac{\sqrt{\tau}}{|\mathcal{S}|}\mathbf{e} = \xi^{\hat\pi}= \tau\cdot\frac{1}{\sqrt{\tau}|\mathcal{S}|}\mathbf{e} \leq \tau \overline \chi^t
	\end{align}
	where we define the auxiliary variable $\overline\chi^t:=\sum_{a \in A}\overline\mu^t(s,a)$ and the last inequality holds from the fact that $\overline\mu^t\in\mathcal{U}$.

	We recall the statement of  \cite[Lemma A.7]{wang2020randomized} and repeat here for quick reference. From the statement of Lemma \ref{intermediate_lemma}, we have
	
	\begin{lemma}[Lemma A.7 \citep{wang2020randomized}]\label{intermediate_lemma}For a given policy $\pi$, the induced stationary distribution $ \xi^{\pi}$, and the average reward $\lambda^\pi$ would satisfy 
		\begin{align}
			\lambda_\pi=\left(\xi^{\pi}\right)^T\sum_{a \in A}\textbf{diag}(\pi_a)\left((P_a-I)v^*+r_a\right)
		\end{align}
		and we have
		\begin{align}\label{lemma_main}
			\lambda^*-\lambda_\pi=\left(\xi^{\pi}\right)^T\sum_{a \in A}\textbf{diag}(\pi_a)\left(\lambda^*\cdot\bbe+(I-P_a)v^*-r_a\right)
		\end{align}
	\end{lemma}
	\normalfont
	Let us first present the proof of Lemma \ref{intermediate_lemma}. Note that since $\xi^\pi$ is the stationary distribution, it holds that $(\xi^\pi)^TP^\pi=(\xi^\pi)^T$. To obtain the first result, we start with the equality that $\lambda^\pi=(\xi^\pi)^Tr$ and then use the fact that $(\xi^\pi)^T(P^\pi-I)=0$ to obtain
	\begin{align}
		\lambda^\pi=&(\xi^\pi)^Tr\nonumber
		\\
		=&(\xi^\pi)^T\left((\xi^\pi)^T(P^\pi-I)v^*+r\right)
		\nonumber
		\\
		=&(\xi^\pi)^T\sum_{a \in A}\textbf{diag}(\pi_a)\left((\xi^\pi)^T(P_a-I)v^*+r_a\right),
	\end{align}
	which holds due to the fact that $(\xi^\pi)^T\sum_{a \in A}\textbf{diag}(\pi_a)\bbe=1$ and implies the statement in \eqref{lemma_main}. Hence, for policy $\hat\pi$, we can write
	%
	%
	%
	\begin{align}
		\lambda^*-\lambda_{\hat\pi}=&\left(\xi^{\hat\pi}\right)^T\sum_{a \in A}\textbf{diag}(\hat\pi_a)\left(\lambda^*\cdot\bbe+(I-P_a)v^*-r_a\right)
		\nonumber
		\\
		=&\frac{1}{T}\sum_{t=1}^{T}\left(\xi^{\hat\pi}\right)^T\sum_{a \in A}\textbf{diag}(\overline\pi^t_a)\left(\lambda^*\cdot\bbe+(I-P_a)v^*-r_a\right)
		\nonumber
		\\
		\leq&\tau\frac{1}{T}\sum_{t=1}^{T}\left(\overline\chi^{t}\right)^T\sum_{a \in A}\textbf{diag}(\overline\pi^t_a)\left(\lambda^*\cdot\bbe+(I-P_a)v^*-r_a\right)
		\nonumber
		\\
	\end{align}
	where we have used the fact that $\xi^{\hat\pi}\leq\tau \overline \chi^t$ and primal feasibility $\left(\lambda^*\cdot\bbe+(I-P_a)v^*-r_a\right)\geq 0$. Next, from the definition of $\overline\pi^t_a$ and $\overline\chi^{t}$, we can write
	\begin{align}
		\lambda^*-\lambda_{\hat\pi}\leq&\tau\left(\frac{1}{T}\sum_{t=1}^{T}\sum_{s \in \mathcal{S}}\sum_{a \in \mathcal{A}}\left(\overline\mu^{t}(s,a)\right)^T\left((v^*-P_av^*)-r_a\right)(s) + \lambda^*\right)
	\end{align}
	Via Markov inequality, we obtain
	\begin{align}
		\lambda^*-\lambda_{\hat\pi}	\leq&\frac{3}{2}\tau\left(\frac{1}{T}\sum_{t=1}^{T}\mathbb{E}\left[\sum_{a \in \mathcal{A}}\left(\overline \mu^{t}(a)\right)^T\left((v^*-P_av^*)-r_a\right)\right] + \lambda^*\right)
	\end{align}
	with probability $2/3$ at least. And further from the statement of Lemma \ref{eq:KL_divergence_martingale_final400}, we can write
	\begin{align}
		\lambda^*-\lambda_{\hat\pi}	\leq&\frac{3}{2}\tau\widetilde{\mathcal{O}}\left( \tilde{t}_{mix}\sqrt{\frac{n\mathcal{E}_{0} |\mathcal{S}||\mathcal{A}|D(\Gamma, \rho)}{T}}\right). 
	\end{align}
	By selecting $T=\Omega\left(\tau^2\tilde{t}_{mix}^2\frac{n\mathcal{E}_{0} |\mathcal{S}||\mathcal{A}|D(\Gamma, \rho)}{\epsilon^2}\right)$, we can write that %
	\begin{align}
		\lambda^*-\epsilon 	\leq\lambda_{\hat\pi}, 
	\end{align}
	with probability $2/3$.
%

\section{Proof of Theorem \ref{eq:KL_divergence_martingale_final400}}\label{theorem_2}

	From Algorithm \ref{alg_single2}, it is clear that each independent $K$ trials generates an $\epsilon/3$ optimal policy. Further, for each $\overline{Y}^{(k)}$,  according to \cite[Lemma A.8]{wang2020randomized}, we can write that for a given policy $\overline{\pi}^{(k)}$, there exists an algorithm such that
	\begin{align}
		\overline{Y}^{(k)}- \lambda_{\overline{\pi}^{(k)}} \in\left[-\frac{\epsilon}{3},\frac{\epsilon}{3}\right],
	\end{align}  	
	holds with probability $1-\frac{\delta}{2K}$ in $L$ number of steps where $L=\tilde{\mathcal{O}}\left(\frac{t_{mix}}{\epsilon^2}\log\left(\frac{4K}{\delta}\right)\right)$. We note that step 3 in Algorithm \ref{alg_single2} takes $KL$  number of steps to execute.  Next, since we select $\widetilde\pi=\overline\pi^{(k^*)}$, we need to show that $\widetilde\pi$ is indeed near optimal with probability $1-\delta$. The number of samples required by the above steps is $\tilde{\mathcal{O}}\left[N_{\epsilon/3}\log\left(\frac{1}{\delta}\right)+L\log\left(\frac{1}{\delta}\right)\right]$ time steps.  
	
	Next, the goal is to prove that the output policy $\hat pi$ of the meta Algorithm \ref{alg_single2} is indeed $\epsilon$ optimal with probability at least $1-\delta$. To achieve that, we need to appropriately select $K$. Let us define $$\mathcal{K}:=\big\{k\in[K]~|~ \lambda_{\overline{\pi}^{(k)}}\geq 	\lambda^* -\frac{\epsilon}{3}\big\},$$
	which denotes the set of all possible successful trials of Algorithm \ref{alg_single2}. Let us consider an event where $\mathcal{K}\neq \emptyset$ and the policy evaluations errors $\overline{Y}^{(k)}- \lambda_{\overline{\pi}^{(k)}}\in \left[-\frac{\epsilon}{3},\frac{\epsilon}{3}\right]$  which implies that we have $$\lambda_{\overline{\pi}^{(k)}} -\frac{\epsilon}{3}\leq \overline{Y}^{(k)} \leq \lambda_{\overline{\pi}^{(k)}}+ \frac{\epsilon}{3},$$
	for all $k$ and we note that $\lambda_{\overline{\pi}^{(k)}}\geq 	\lambda^* -\frac{\epsilon}{3}$ if $k\in\mathcal{K}$. For $\mathcal{K}$ being non-empty, the output policy with the largest value $\overline{Y}^{(k)} $ value would satisfies $$\overline{Y}^{(k)} \geq \lambda^* -\frac{2\epsilon}{3}.$$ Since the policy evaluation error is upper bounded by $\epsilon/3$, this would imply that the policy $\widetilde\pi$ is $\epsilon$ optimal. Let us consider the probability of event $\left(\lambda_{\widetilde\pi} \leq  \lambda^*-{\epsilon}\right)$ as 
	\begin{align}
		\mathbb{P}\left(\lambda_{\widetilde\pi} <  \lambda^*-{\epsilon}\right) \leq &	\mathbb{P}\left(\{\mathcal{K}=\emptyset\} ~ \cup ~ \{\exists k : 	\overline{Y}^{(k)}- \lambda_{\overline{\pi}^{(k)}} \notin\left[-\frac{\epsilon}{3},\frac{\epsilon}{3}\right]\}\right)
		\nonumber
		\\
		\leq &	\mathbb{P}\left(\mathcal{K}=\emptyset \right)+  \left(\exists k : 	\overline{Y}^{(k)}- \lambda_{\overline{\pi}^{(k)}} \notin\left[-\frac{\epsilon}{3},\frac{\epsilon}{3}\right]\right).
	\end{align}
	From the definition of $\mathcal{K}$, we can write
	\begin{align}
		\mathbb{P}\left(\lambda_{\widetilde\pi} <  \lambda^*-{\epsilon}\right)
		\leq &\prod_{k=1}^{K} \mathbb{P}\left(\lambda_{\overline{\pi}^{(k)}}< 	\lambda^* -\frac{\epsilon}{3}\right)+  \sum_{k=1}^{K} \mathbb{P}\left(\overline{Y}^{(k)}- \lambda_{\overline{\pi}^{(k)}} \notin\left[-\frac{\epsilon}{3},\frac{\epsilon}{3}\right]\right)
		\nonumber
		\\
		\leq \left(\frac{1}{3}\right)^K + \frac{\delta}{2}.
	\end{align}
	By selecting $K=\log_{1/3} \left(\frac{\delta}{2}\right)$, we can write that 
	\begin{align}
		\mathbb{P}\left(\lambda_{\widetilde\pi} <  \lambda^*-{\epsilon}\right)
		\leq \delta.
	\end{align}
	Hence proved. 

\section{Additional Experiments}
\begin{figure*}[htb!]
	\centering
	\subfigure[  Occupancy measure consensus error.]{
			\includegraphics[scale = 0.20]{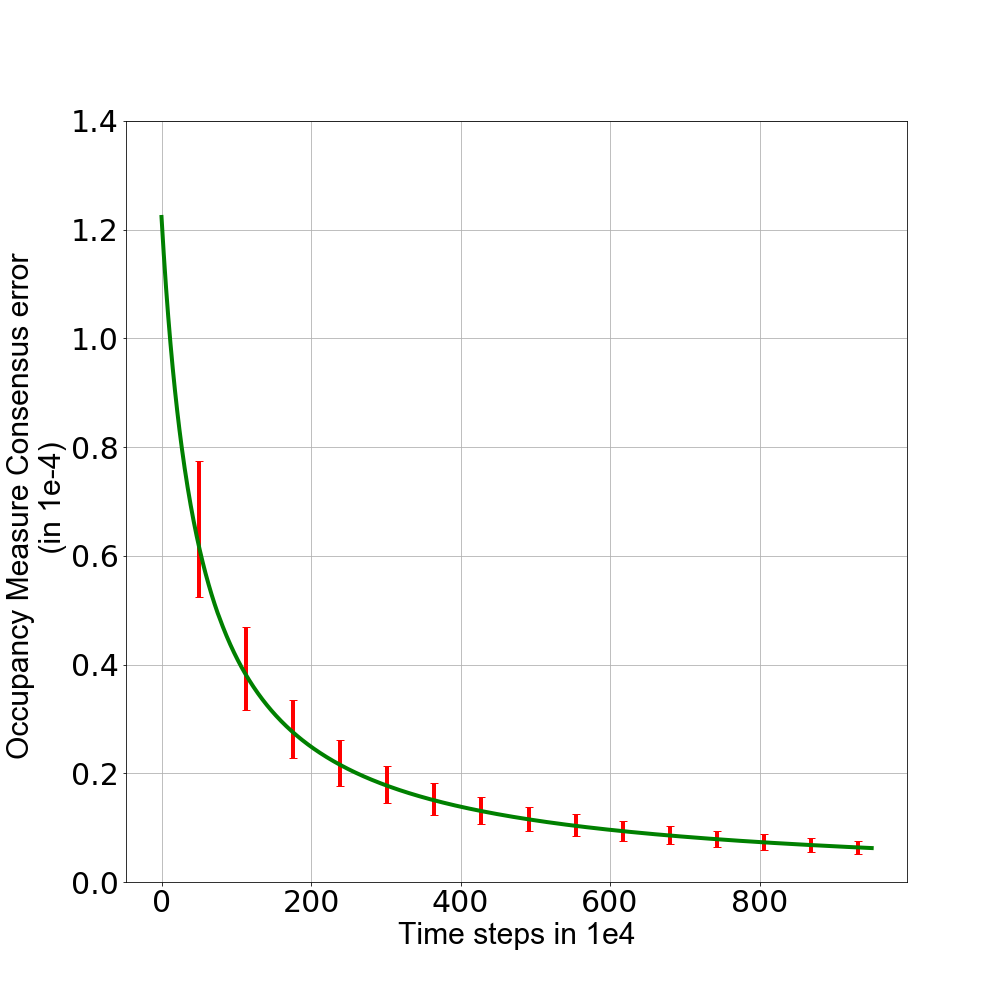}	
		\label{fig:regret_rs_1}
	}
		\hspace{.25in}
\subfigure[  Value function consensus error.]{
\includegraphics[scale = 0.20]{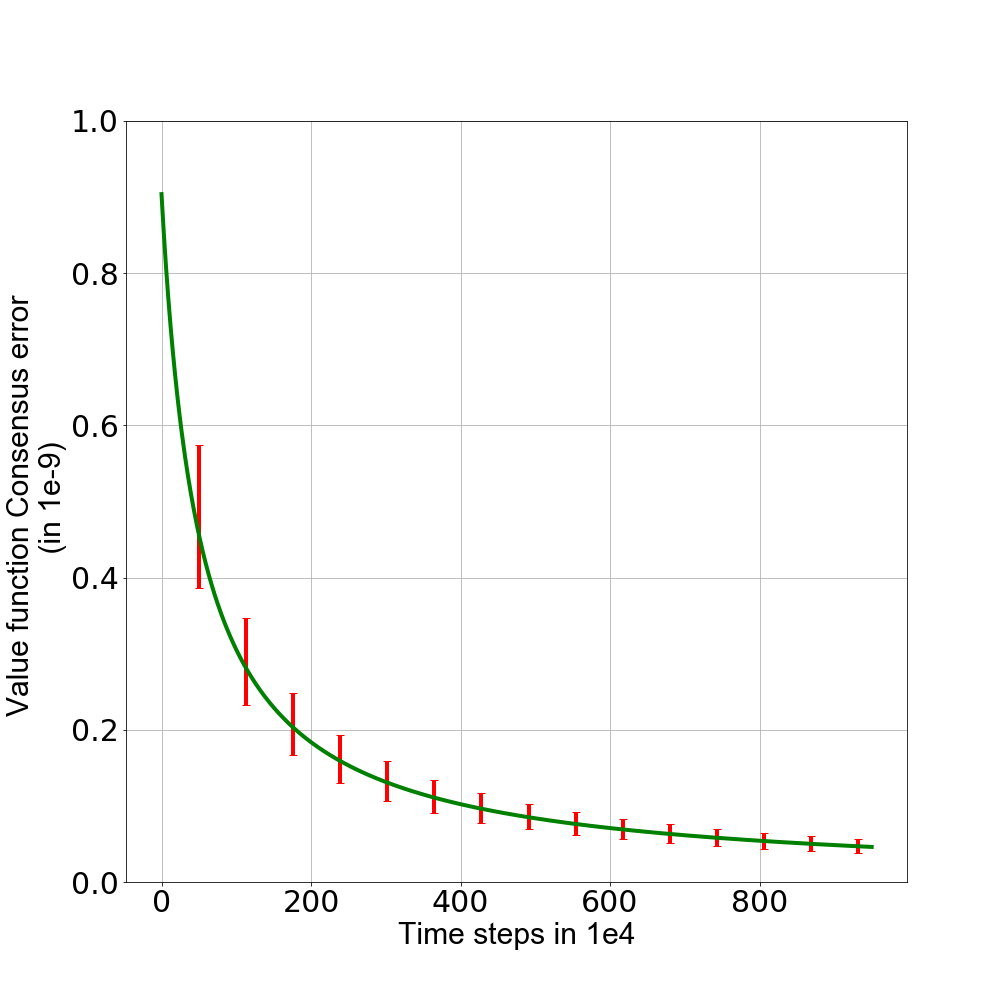}
		\label{fig:regret_2}
	}    
	\caption{Consensus error plots for primal-dual pairs $(\widetilde{\mu}_i^t,\widetilde{v}_i^t)$. Fig. \ref{fig:regret_rs_1} plots consensus error for occupancy measure. Fig. \ref{fig:regret_2} plots consensus error for value function.   }
	\label{fig:2}
\end{figure*}

In this section, we plot the consensus error in Fig. \ref{fig:2} for the primal-dual pair $(\widetilde{\mu}_i^t,\widetilde{v}_i^t)$ for the multi-agent experimental setting explained in Sec. \ref{sec_experiments}.  As mentioned in \eqref{eq:consensus_round}, each agent $i$ shares its primal-dual variables with others according to network graph $\mathcal{G}$, and in turn performs weighted averaging to update its local estimates. In our experiments, we use the Erdos-Reyni model to generate random network graph $\mathcal{G}$ among the agents to perform consensus. Specifically, with $n=3$ agents we utilized binomial model to generate graphs with edge probability of $p=0.3$.

\end{document}